%% file: main.tex
\documentclass[hidelinks,onefignum,onetabnum]{siamonline220329}

\usepackage{lipsum}
\usepackage{amsmath,amsbsy,amsfonts,amssymb}
\usepackage{graphicx}
\usepackage{xcolor}
\usepackage{mathtools}
\usepackage{epstopdf}
\usepackage{algorithmic}
\usepackage{soul}
\ifpdf
  \DeclareGraphicsExtensions{.eps,.pdf,.png,.jpg}
\else
  \DeclareGraphicsExtensions{.eps}
  \fi

\usepackage{tikz}
\usetikzlibrary{shapes,arrows}

\input{defs.tex}
\input{trimFig.tex}
\usepackage{empheq}
\usepackage{placeins}
\usepackage{hyperref}
\hypersetup{bookmarks=true,colorlinks=true, linkcolor=blue,urlcolor=red,citecolor=siaminlinkcolor,linktoc=all}
\usepackage{listings}
\usepackage{array}
\usepackage[noadjust,compress]{cite}

\newcommand{\xt}{\tilde{x}}
\newcommand{\yt}{\tilde{y}}
\newcommand{\yh}{\hat{y}}
\newcommand{\Bt}{\tilde{B}}
\newcommand{\Ct}{\tilde{C}}
\newcommand{\gt}{\tilde{g}}

\newcommand{\tslow}{\tau}

\newcommand{\remove}[1]{}

\headers{Fast-slow neural networks}{D. A. Serino, A. Alvarez Loya, et al.}

\title{{Fast-Slow Neural Networks for Learning Singularly Perturbed Dynamical Systems}\thanks{Submitted to the editors February 2024.
Corresponding author: \email{qtang@lanl.gov}
\funding{
 This work was partially supported by the U.S. Department of Energy Advanced Scientific Computing Research (ASCR) under DOE-FOA-2493 ``Data-intensive scientific machine learning''. 
 It was also partially supported by the ASCR program of Mathematical Multifaceted Integrated Capability Center
(MMICC). 
Alvarez Loya is supported by the U.S. National Science Foundation under NSF-DMS-2213261.
}}}
\author{Daniel A. Serino\thanks{Los Alamos National Laboratory, Los Alamos, NM 
  (\email{dserino@lanl.gov},
   \email{aalvarezloya@lanl.gov}).}
\and Allen Alvarez Loya\footnotemark[2]
\and J. W. Burby\thanks{Los Alamos National Laboratory, Los Alamos, NM. Current address: Department of Physics, University of Texas at Austin, Austin, TX (\email{joshua.burby@austin.utexas.edu}).}
\and Ioannis G. Kevrekidis\thanks{Department of Chemical and Biomolecular Engineering and Department of Applied Mathematics and Statistics, Johns Hopkins University, Baltimore, MD (\email{yannisk@jhu.edu}).
}
\and Qi Tang\thanks{Los Alamos National Laboratory, Los Alamos NM. Current address: School of Computational Science and Engineering, Georgia Institute of Technology, Atlanta, GA. (\email{qtang@lanl.gov}, \email{qtang@gatech.edu})}
}

\usepackage{amsopn}

\ifpdf
\hypersetup{
  pdftitle={FSNN},
  pdfauthor={Daniel A. Serino et. al.}
}
\fi

\newcommand{\bogus}[1]{}

\begin{document}

\maketitle

\begin{abstract}
Singularly perturbed dynamical systems 
play a crucial role in climate dynamics and 
plasma physics. 
A powerful and well-known tool to address these systems is the Fenichel normal form, which significantly simplifies fast dynamics near slow manifolds through a transformation. However, this normal form is difficult to realize in conventional numerical algorithms. In this work, we explore an alternative way of realizing it through structure-preserving machine learning.
Specifically, a fast-slow neural network (FSNN) is proposed for learning data-driven models of singularly perturbed dynamical systems with dissipative fast timescale dynamics. Our method enforces the existence of a trainable, attracting invariant slow manifold as a hard constraint.
Closed-form representation of the slow manifold 
enables efficient integration on the slow time scale
and significantly improves prediction accuracy beyond the training data.
We demonstrate the FSNN on examples including the Grad moment system, two-scale Lorenz96 equations, and Abraham-Lorentz dynamics modeling radiation reaction of electrons.

\end{abstract}

\begin{keywords}
Singularly perturbed dynamical systems, slow manifold, model reduction, closures
\end{keywords}



\section{Introduction}
\label{sec:Intro}

In this article, a new architecture called the Fast-Slow Neural Network (FSNN) is developed for learning data-driven models of singularly perturbed dynamical systems with dissipative fast timescale dynamics. 
The constraint is imposed by way of a novel neural network architecture with structured weight matrices. By leveraging the neural ODE framework from~\cite{chen2018}, our technique produces deterministic evolution equations that approximately reproduce trajectory training data. When the training data comes from an unknown dynamical system, the method achieves data-driven \emph{full-order} model discovery with stability guarantees due to provable existence of the slow manifold. 
In short, our architecture is, by design, a singularly perturbed dynamical system which contains an attractor.

Various dynamical systems encountered in nature exhibit a timescale separation such that the fast timescale is dissipative. Prominent examples include the kinetic dynamics of strongly-collisional gasses and plasmas, and radiative transfer in optically-thick media. For ordinary differential equations (ODEs) of this type, Fenichel theory~\cite{fenichel} reveals that fast timescale dissipation implies existence of an invariant slow manifold that attracts nearby trajectories and governs the long-term dynamics. For partial differential equations (PDEs), or more generally integro-differential equations, Fenichel theory generally only applies formally, but its predictions are often valid anyway. 
For PDEs, inertial manifold theory~\cite{foias1988inertial} has been studied and shown to be useful in applications such as model reductions of fluid dynamics.

Nonlinear dimensionality reduction
inspired by inertial manifold 
theory has a rich history (see~\cite{koronaki2023nonlinear} and the references therein)
originating in the 1980s and 1990s~\cite{foias1998a,foias1988inertial,foias1989,jolly1989,temam1989a,temam1989b,jauber1990,guermond08}.
Data-driven approaches for discovering inertial manifolds in dynamical systems 
and using separations of scales to evolve a lower-dimensional reduced order model first materialized in the 90s~\cite{jolly1990approximate,titi1990,kramer1991,krischer93,shvarts1998,theo2000,yannis2000}
and recently reemerged in popularity due to boom in machine-learning-based 
modeling~\cite{benner2015,chorin2015,lu2017,sinai2019,lee2020,linot2020,zeng2022,linot2022,zeng2023autoencoders,jesus2023,linot2023couette,linot2023, patsatzis2023slow}.
The work in this manuscript tackles a special class of systems with attracting invariant sets, namely fast-slow systems with dissipative fast timescale dynamics~\cite{burby2020}.
Unlike previous data-driven 
approaches, the FSNN, {a special form of neural ODE \cite{chen2018}}, is
a \emph{full-order} model
which, using 
neural network architectures with special properties,
strictly enforces the existence of an attracting slow manifold.
This method may be viewed as a major step along a pathway in developing data-driven model discovery techniques with embedded invariant manifolds.
{In future work, we plan to address the general case where 
explicit timescale separation is not used as a crutch.}

The secondary goal of the proposed network is to discover an efficient and asymptotic-preserving \emph{reduced-order} model.
In the common methodology of physics-based model reduction, one first decides a latent space  (e.g. lower-order moments of a distribution) and then looks for a closure using tools such as asymptotic analysis. 
Many complicated closures have been discovered through this approach. 
The classical examples in rarefied gases and plasma physics include Grad's moment method~\cite{grad13kinetic}, Braginskii closure~\cite{braginskii1965transport} and Hammett-Perkins closure~\cite{hammett1990fluid}.
The challenge of this methodology is that an exact closure may not exist. The accuracy of the approximate closure 
is related to the dimension of the latent space.
In contrast, Fenichel theory guarantees the existence
of an exact closure in the neighborhood of the 
singular limit of the equations
in the form of a diffeomorphic coordinate transformation.
This transformation exposes an invariant manifold
by flattening the fast variable space.
Often in practice, as asymptotic expansion is the best
analytical approximation to the exact closure.
In our work, we leverage structured neural network
architectures and optimization to learn an approximation 
to the exact closure that is hidden in the dynamics.
The data-driven approach described in this work can discover such a transformation only using short-term trajectory data.
Therefore, this approach may be useful when the underlying dynamical system is known but the targeting latent space or a formal closure are not known.

The design of the FSNN architecture is motivated by 
the first few theorems discovered by Fenichel~\cite{Jones1995GeometricSP}.
The components of the architecture are constructed 
using structured weight matrices.
The coordinate transformation into Fenichel normal form is based on an invertible coupling flow network~\cite{teshima20} 
in which we introduced the additional structure-preserving property of bi-Lipschitz control through the use of
bi-Lipschitz affine transformation (bLAT)
layers. 
This network and its inverse can be evaluated using closed-form expressions with identical computational complexity.
We note that ``vanilla''
neural networks are generally 
non-invertible~\cite{noninv1998,noninv2000,noninv2023}.
We also introduce the Schur form network, which is a novel 
architecture based on the Schur decomposition for controlling the locations of eigenvalues 
in the complex plane, to parameterize the 
linear term in the dynamics.
Additionally, a low-rank {bilinear map} network is included to increase the approximation power in the fast variable space.
{Our approach is then integrated 
using adaptive time stepping
that treats the linear term implicitly,
further leveraging back-substitution on the Schur form.}
Some of the structures introduced are related to the idea of trivializing manifold optimization through parameterizations~\cite{lezcano2019trivializations}.
Structure-preserving flow maps through machine learning 
are an active area of research, and include
such as the symplectic neural networks 
for Hamiltonian dynamics~\cite{sympnet2020, burby2021}. 
Previous studies on flow map learning for multi-scale dynamical systems include~\cite{liu2022hierarchical} and~\cite{duruisseaux2023approximation}.
In particular,~\cite{duruisseaux2023approximation} simultaneously preserves symplecticity and adiabatic invariance in a multi-scale setting.


Before going into a detailed presentation of our new method, it is worth highlighting an important shortcoming it suffers from, and how that shortcoming might be mitigated in practice. 
The method assumes that the dimension of the attracting invariant manifold is known in advance. While this is not a limitation in the surrogate modeling context, where the (expensive to simulate) dynamical system is known in detail, it does present a moderate challenge in the model discovery context. As such, when performing data-driven model discovery using our method, it will be necessary to first determine the dimension of the slow manifold, either through physical reasoning, or through application of auxiliary data mining techniques
such as diffusion maps~\cite{coifman2006diffusion},
autoencoders~\cite{zeng2022}, or other 
techniques from reduced order modeling~\cite{koronaki2023nonlinear}. 
For a review on various manifold learning techniques, we refer the reader to~\cite{manifoldlearning}.
In many cases, including the examples considered in this manuscript, 
the form of the slow manifold and therefore its dimensionality are known in the singular 
limit $\epsilon=0$ as a power series expansion of $\epsilon$. 
However, the exact form for $\epsilon>0$ and the transformation into Fenichel normal form 
is unknown.
In order to learn the slow manifold effectively, a sufficient number of trajectories
for various choices of $\epsilon$ and initial conditions are necessary.
In the final example of the paper, we apply our network in the nontrivial setting
of learning the slow manifold for the Abraham-Lorentz equations which,
for trajectories starting off of the slow manifold, are unstable forward in time.

The remainder of the manuscript is organized as 
follows.
A review of fast-slow dynamical systems is presented in
Section~\ref{sec:fastslowfenichel}.
The fast-slow neural network architecture is 
described in detail in Section~\ref{sec:fsnn}.
The numerical implementation of the architecture
is presented in Section~\ref{sec:numerical}.
Results from applying the architecture 
to various examples are shown in Section~\ref{sec:examples}.
Concluding remarks are made in Section~\ref{sec:conclusions}.


\section{Fast-Slow Dynamics and the Fenichel Normal Form}
\label{sec:fastslowfenichel}

In this paper, we develop a novel machine learning (ML)-inspired
architecture, the Fast-Slow Neural Network (FSNN). 
FSNNs are a parameterization for 
a class of dynamical systems exhibiting multiscale behavior in time
that, to leading order, are dissipative on the fast time scale.
Fast-slow systems, which exhibit both short-term and long-term 
behavior, are relevant to many physical systems observed in nature,
including the climate and plasma systems.
In this section, we motivate the FSNN architecture
with a review of the theory of fast-slow systems.

{\definition
A fast-slow system is an ODE of the form~\cite{Jones1995GeometricSP, burby2020}~
\bse
\label{eq:fastslow}
\begin{align}
\epsilon \frac{d}{d\tslow} y &= f(x, y, \epsilon), \\
\frac{d}{d\tslow} x &= g(x, y, \epsilon),
\end{align}
\ese
where $x=x(\tslow)\in\mathbb{R}^{N_x}$, $y=y(\tslow)\in\mathbb{R}^{N_y}$, 
$f$ and $g$ are smooth functions of their arguments,
and $D_y f(x, y, 0)$ is invertible for all $(x, y)$ where $f(x, y, 0) = 0$. \label{def:fastslow}
}

\medskip

\noindent In the above definition, 
$y$ and $x$ are the fast and slow variables, respectively, 
over corresponding timescale variables $t$ and $\tau=\epsilon t$.

In some instances, a dynamical system may exhibit 
a fast-slow split, but a representation for the fast and slow variables 
is unknown.
Consider the following representation of a 
dynamical system for $z(\tau)\in\mathbb{R}^N$,
$\epsilon>0$,
\begin{align}
\epsilon \frac{d}{d\tau} z &= U(z, \epsilon),
\label{eq:original}
\end{align}
where $U$ is a sufficiently smooth function of $z$ and $\epsilon$.

{\definition
The dynamical system~\eqref{eq:original} exhibits a fast-slow split if
there is an $\epsilon$-de\-pendent invertible transformation
of the dependent variables 
into fast-slow coordinates $(x, y)$, $x\in\mathbb{R}^{N_x}$, 
$y\in\mathbb{R}^{N_y}$, $N = N_x + N_y$, 
where the dynamics are described by a fast-slow system.
\label{def:invertibletrans}
}

\medskip

\noindent In this paper, we restrict our attention to the
\textit{normally stable fast-slow system},
which is a fast-slow system 
where {the fast directions are stable.}

{\definition
A normally stable fast-slow system is a fast-slow system 
where the eigenvalues of the Jacobian, 
$D_y f(x, y, 0) \in \mathbb{R}^{N_y\times N_y}$,
are strictly to the left of the imaginary axis when $f(x, y, 0) = 0$.
\label{def:normfastslow}
}

\medskip

\noindent Due to the dissipative fast-scale dynamics
of normally stable fast-slow systems,
the solutions will be attracted to an
invariant slow manifold in forward time.
We will now define the notion of the slow manifold for 
normally stable fast-slow systems.
Consider a fast-slow system in the form~\eqref{def:fastslow}.
{\lemma For $\epsilon>0$ sufficiently small, there exists a 
function, $y^*(x, \epsilon)$, such that the graph
\begin{align}
M_\epsilon = \{(x, y): y = y^*(x, \epsilon)\},
\end{align}
is locally invariant under~\eqref{def:fastslow}. 
We define $M_\epsilon$ to be the slow manifold.
\label{lem:slowmanifold}
}

\medskip

\noindent Here a bounded hypersurface is called a \emph{locally invariant} manifold
 if the vector field defining the dynamical system is tangent to the hypersurface at all of its points \cite{sawant2006model}.
This lemma is proved following 
Theorem 4 in~\cite{Jones1995GeometricSP}.
The purpose of the FSNN is to be a
universal approximator for normally stable fast-slow systems. 
The network is based on the following
well known result due to Fenichel.

\medskip

{\lemma
When $\epsilon>0$ is sufficiently small and $(x, y)$ lies
is in a neighborhood around the slow manifold $M_\epsilon$,
there exists a change of coordinates
$(x, y) \rightarrow (\xt, \yt)$
where the dynamics of normally stable fast-slow systems 
can be written as,
\bse
\label{eq:fenichelnormal}
\begin{align}
\frac{d}{d t} \yt &= A(\xt, \yt, \epsilon) \yt, \\
\frac{d}{d t} \xt &= \epsilon {\gt}(\xt, \yt, \epsilon),
\end{align}
\ese
where $A(\xt, \yt, \epsilon)\in\mathbb{R}^{N_y\times N_y}$,
${\gt}(\xt, \yt, \epsilon)\in\mathbb{R}^{N_x}$,
and $A(\xt, 0, 0)$ is a stable matrix, i.e., its eigenvalues are on the left of the imaginary axis.
}

\medskip

\noindent The above lemma is proved in Chapter~3.2
of~\cite{Jones1995GeometricSP}.
The lemma  introduces \eqref{eq:fenichelnormal}, which is often called the \emph{Fenichel normal form} due to the pioneering work of Neil Fenichel~\cite{fenichel}.
\label{lem:fenichelnormal}
Note that the Fenichel normal form given in~\cite{Jones1995GeometricSP} is more general than \eqref{eq:fenichelnormal}, {including
fixed points that are either stable or unstable in the linearization of the fast variable.}
In this work, we only consider fast variables that are long-time stable. 
Moreover, due to the existence of the coordinate transformation, 
we consider the general form of the singularly perturbed systems~\eqref{eq:original}. 
However, the advantage of the split form~\eqref{eq:fastslow} is 
that the slow manifold dimension is known. 

Our FSNN takes advantage of a modified formulation
of Fenichel normal form which will be described in the 
following discussion.
The following theorem forms the framework for the 
FSNN architecture.

{\theorem
Suppose the dynamics of $z(t;\epsilon)\in\mathbb{R}^N$, $\epsilon>0$,
admit a fast-slow split 
and can be described by normally stable fast-slow system.
Then there exists a $\epsilon$-dependent invertible transformation
\begin{align}
    \left[\begin{array}{c}
                y \\ x
      \end{array}\right] &= h(z; \epsilon),
      \label{eq:transformation}
\end{align}
where $x(t) \in \mathbb{R}^{N_x}$, $y(t)\in \mathbb{R}^{N_y}$, 
are slow and fast variables, respectively, 
and for $|y| \le \Delta$, where $\Delta>0$,
the dynamics for $(x, y)$ can be described using
\begin{subequations}
\begin{align}
    \frac{d}{dt} y &= T(x) y + B(x, y, \epsilon)(y, y)
    + \epsilon C(x, y, \epsilon) y, \label{eq:dydt} \\
    \frac{d}{dt} x &= \epsilon g(x, y, \epsilon).
\end{align}
\label{eq:fastslowthm}
\end{subequations}
$T(x)\in\mathbb{R}^{N_y\times N_y}$ 
is a block upper triangular matrix consisting of  
either $\mathbb{R}^{2\times 2}$ or $\mathbb{R}$ blocks on the diagonal
and the eigenvalues of these diagonal blocks are 
strictly to the left of the imaginary axis for all $x$.
$B(x, y, \epsilon): \mathbb{R}^{N_y} \times \mathbb{R}^{N_y} \rightarrow \mathbb{R}^N_y$ is a bilinear map,
$C(x, y, \epsilon)\in\mathbb{R}^{N_y\times N_y}$ is a matrix,
and $g(x, y, \epsilon) \in \mathbb{R}^{N_x}$.
$T$, $B$, $C,$ and $g$ are $C^r$ for $r\ge1$.
\label{thm:fsnn}
}

\medskip

\noindent In the following discussion, 
we will prove Theorem~\ref{thm:fsnn}
by refining the Fenichel normal form in~\eqref{eq:fenichelnormal}
and applying the matrix Schur decomposition.
%
First, we modify~\eqref{eq:fenichelnormal} by performing a Taylor 
expansion of the term $A(\xt, \yt, \epsilon) \yt$
about $\yt=0$ and $\epsilon=0$. 
{Assuming that $A$ is in $C^r$ for $r\ge 1$}, we have
\begin{align}
    A(\xt, \yt, \epsilon) \yt = 
    A(\xt, 0, 0) \yt
    + B(\xt, \yt, \epsilon) (\yt, \yt)
    + \epsilon C(\xt, \yt, \epsilon) \yt ,
    \label{eq:taylorexpansion}
\end{align}
where $B(\xt, \yt, \epsilon):\mathbb{R}^{N_y} \times \mathbb{R}^{N_y} \rightarrow \mathbb{R}^{N_y}$ is a bilinear map 
and $C(\xt, \yt, \epsilon)\in\mathbb{R}^{N_y \times N_y}$
are remainder terms.
Next, the Schur decomposition of $A(\xt, 0, 0)$
will be used to define a new set of fast coordinates.

{\definition The Schur form of a matrix, 
$A\in\mathbb{R}^{M\times M}$,
is a matrix $T\in\mathbb{R}^{M\times M}$
that {forms the Schur decomposition}
$A = Q T Q^T$, where $Q\in\mathbb{R}^{M\times M}$ is
orthogonal ($Q^T Q = I$) and 
$T$ is block upper triangular where
the diagonal blocks consist of $\mathbb{R}$
and $\mathbb{R}^{2\times 2}$ blocks.
}

\bogus{
{\lemma The Schur decomposition of a matrix
$A\in\mathbb{R}^{M\times M}$, is given by
\begin{align}
A = Q T Q^{\rm T}, \label{eq:schur}
\end{align}
where $Q, T\in\mathbb{R}^{M\times M}$, $Q$ is an orthogonal matrix 
($Q^{\rm T} Q = I$), and $T$ is the Schur form of $A$.
The eigenvalues of $A$
are given by the eigenvalues of the diagonal blocks of $T$.
}}

\medskip

\noindent The existence of the 
Schur decomposition is proven in many references, including~\cite{horn_johnson_1985,matrixanalysis},
though for any given matrix, the decomposition is 
not unique.
{However, the eigenvalues of $A$
are given by the eigenvalues of the diagonal blocks of $T$.}
Let the Schur decomposition of $A(\xt, 0, 0)$ be given by
\begin{align}
A(\xt, 0, 0) = Q(\xt) T(\xt) Q(\xt)^{\rm T},
\end{align}
where $Q(\xt)\in\mathbb{R}^{N_y\times N_y}$ is an orthogonal matrix 
and $T(\xt)\in\mathbb{R}^{N_y\times N_y}$ is the Schur form of $A(\xt, 0, 0)$.
Due to our assumptions on $A(\xt, 0, 0)$, all of the eigenvalues of $T(\xt)$ are
to the left of the imaginary axis.
Consider the coordinate transformation $\yh = Q(\xt)^{\rm T} \yt$.
We form the relevant ODE for $\yh$ to use in place of the ODE for $\yt$.
\begin{align*}
\frac{d}{dt}\yh =\;& 
Q(\xt)^{\rm T} \frac{d}{dt} \yt
+ \left( \grad Q(\xt)^{\rm T}\right) \left(\frac{d}{dt} \xt, \yt \right) \\
=\;&{Q(\xt)^T } Q(\xt) T(\xt) Q(\xt)^T \yt 
+ Q(\xt)^{\rm T} B(\xt, Q(\xt)\yh, \epsilon)(Q(\xt)\yh, Q(\xt)\yh) \\
&+ \epsilon Q(\xt)^{\rm T} C(\xt, \epsilon) Q(\xt)\yh
+ \epsilon \left(\grad Q(\xt)^{\rm T}\right) \left(g_\epsilon(\xt, \yh), \yh\right). \\
=\;& T(\xt) \yh 
+ Q(\xt)^{\rm T} B(\xt, Q(\xt)\yh, \epsilon)(Q(\xt)\yh, Q(\xt)\yh) \\
&+ \epsilon Q(\xt)^{\rm T} C(\xt, \epsilon) Q(\xt)\yh
+ \epsilon \left(\grad Q(\xt)^{\rm T}\right) \left(g_\epsilon(\xt, \yh), \yh\right).
\end{align*}
By redefining terms, it can be shown that the ODE is equivalent to 
\begin{align*}
    \frac{d}{dt}\yh =
  T(\xt) \yh 
+ \Bt(\xt, \yh, \epsilon)(\yh, \yh) 
+ \epsilon \Ct(\xt, \yh, \epsilon) \yh,
\end{align*}
where $\Bt(\xt, \yh, \epsilon)$ is a {bilinear map} and
$\Ct(\xt, \yh, \epsilon)$ is a {linear map}.
This establishes an equivalence to~\eqref{eq:fastslowthm}.

Theorem~\ref{thm:fsnn} takes advantage of many useful 
properties. In the fast-slow coordinate system of
Theorem~\ref{thm:fsnn}, the definition of the slow-manifold 
simplifies significantly.
{\theorem The graph,
$M  = \{(x, y): x\in\mathbb{R}^{N_x},\, y=0\}$
is locally invariant under~\eqref{eq:fastslowthm}. 
}

\medskip

\noindent This theorem can be easily verified by substituting
$y=0$ into~\eqref{eq:dydt} and observing that
$\frac{d}{dt} y = 0$ on the slow manifold.
Trajectories on the slow manifold will remain on the slow manifold
for all time. 
In this case, the dynamics  can be reduced to an ODE for the slow variable,
\begin{align}
    \frac{d}{d\tau} x = g(x, 0, \epsilon), 
    \qquad (x, y) \in M,
    \label{eq:slowmanifoldinteg}
\end{align}
where $\tau$ is the slow timescale.
In the case when the coordinates $(x, y)$ are close to, 
but not on $M$,
the following theorem states that trajectories 
emanating from $(x, y)$ will approach $M$ exponentially in 
forward time.
{\theorem If $\epsilon>0$ is sufficiently small, then, for 
$|y| \le \Delta$, where $\Delta > 0$, 
there are $\kappa>0$ and $\alpha<0$ such that 
\begin{align}
    d\left((x(t), y(t)), M\right) \le \kappa \exp(\alpha t),
\end{align}
where $d$ represents the Euclidean distance.
}

\medskip

\noindent This theorem is proved following Theorem 5 in~\cite{Jones1995GeometricSP}.

Theorem~\ref{thm:fsnn} provides a framework
for the FSNN architecture,
which expresses each of the functions 
$h, T, B, C,$ and $g$ using 
neural networks with specially designed architectures
described in the following section.
There are a few considerations that went into using the
modified Fenichel normal form. 
By expanding $A(\xt, \yt, \epsilon)$
in~\eqref{eq:taylorexpansion}, we expose three separate
terms with known structures.
We further use the Schur decomposition to 
explicitly expose the eigenvalues of $A(\xt, 0, 0)$,
which are required to be left of the imaginary axis.
Finally, we further simplify the linear term by absorbing 
the matrix $Q(\xt)$ into the representation.
Though we use the Schur decomposition
to represent $A(\xt, 0, 0)$ in this paper, 
other decompositions
for representing matrices with eigenvalues to the 
left of the imaginary axis have been studied in 
literature~\cite{stablematrix, matrixanalysis}.


\section{Fast-Slow Neural Networks (FSNNs)}
\label{sec:fsnn}

\begin{figure}
    \centering
    \includegraphics[width=1\textwidth]{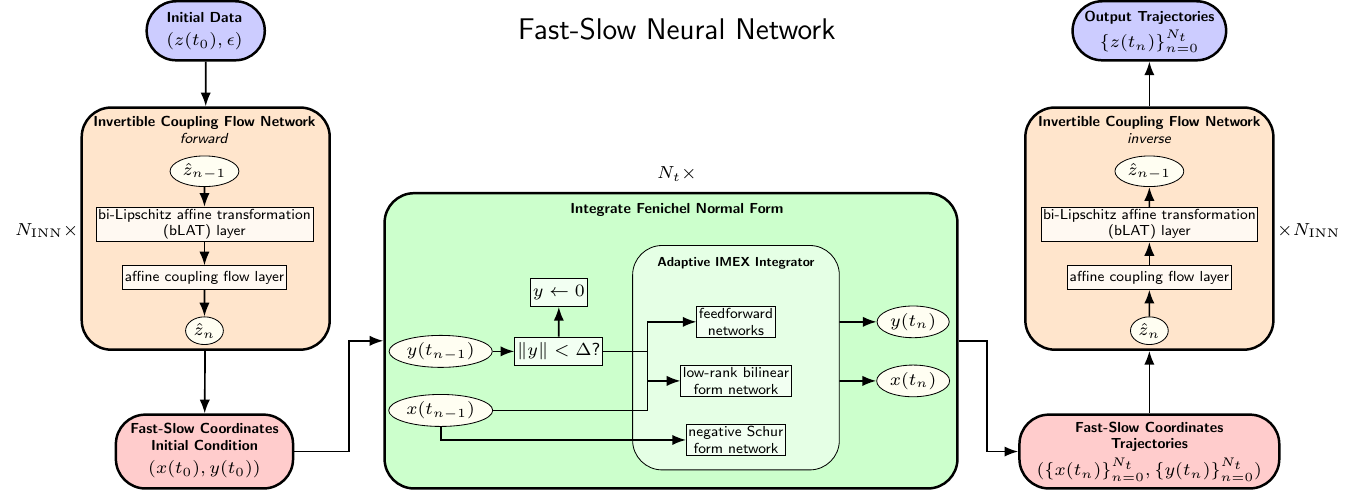}
    \caption{Schematic of the fast-slow neural network (FSNN).
    The original coordinates are projected into fast-slow 
    coordinates using the invertible coupling flow network. 
    In these coordinates, numerical integration 
    of the Fenichel normal form~\eqref{eq:fastslowthm} is used 
    to advance the solution to later time steps. 
    The inverse function of the invertible coupling flow network is used to convert the trajectories in the
    fast-slow coordinates into the original coordinate system.
    }
    \label{fig:FSNN}
\end{figure}

In this section, we present a neural network-based 
approach to approximating the set of solutions to
normally stable fast-slow systems using the framework of
Theorem~\ref{thm:fsnn}.
We refer to this approach as the fast-slow neural network (FSNN),
which is defined below and summarized in Figure~\ref{fig:FSNN}.

{\definition
A fast slow neural network (FSNN)  is a special neural ODE that represents a pullback vector field on a phase space $Z\ni z$ of the form $V(z) = (Dh(z))^{-1}v(h(z))$, where $h$ is a parametric diffeomorphism and $v$ is a parametric vector field in the modified Fenichel normal form \eqref{eq:fastslowthm}. The quantities $h,T,C,B,g$ are parameterized as follows. The diffeomorphism $h$ is represented as an invertible coupling flow network; $T$ is represented using a negative Schur form network; $B$ is represented using a low-rank bilinear map network; and $C,g$ are represented using multi-layer perceptrons.\label{def:fsnn}
}

\noindent{Note that the corresponding flow maps $F_{\Delta t}$ for $V$ and $f_{\Delta t}$ for $v$ are related according to $F_{\Delta t} = h^{-1}\circ f_{\Delta t}\circ h$. It follows that computing the flow map $F_{\Delta t}$ from a FSNN reduces to first computing the flow map $f_{\Delta t}$ for $v$ and then appropriately applying $h$ or $h^{-1}$.} The key components of Definition~\ref{def:fsnn}, including the
invertible coupling flow network,
Schur form network,
and low-rank {bilinear map} network, will be defined in the following
sections.
First we will begin by stating the main result of the 
paper.



{\theorem The FSNN is a universal approximator for 
the set of solutions of normally stable fast-slow systems in the neighborhood of the slow manifold.
\label{thm:fsnnapprox}
}

\medskip

\noindent Due to Theorem~\ref{thm:fsnn}, 
normally stable fast-slow systems can be represented 
using~\eqref{eq:transformation}--\eqref{eq:fastslowthm}. 
Therefore, the FSNN inherits the universal approximation properties
of its components, which we will later show are 
universal approximators for their corresponding class of 
functions. 
The invertible coupling flow network which represents
$h$ is described in Section~\ref{sec:inn} and was proven to be
a universal approximator to invertible maps 
in~\cite{teshima20}.
The universal approximation theorems of the negative Schur form network representing $T$ 
and the {bilinear map} network representing $B$
are proved in Sections~\ref{sec:schur} and~\ref{sec:bilinear},
respectively.
The universal approximation of nonlinear functions
using feedforward neural networks for $C$ and $g$ are 
well known, e.g.~\cite{approximation}.
Finally, we describe a consistent numerical scheme to obtain 
trajectories in time, which has an error bound that 
vanishes as the time steps $\Delta t_n = t_{n}-t_{n-1} \rightarrow 0$.

Our ML-based approach naturally comes with a set of hyperparameters. 
The standard hyperparmeters are those related to the architectures of $h$, $T$, $C$ and $g$.
The only systemic hyperparameters are the dimensions of 
the slow and fast subspaces, $N_x$ and $N_y$, respectively. 
In the examples considered in Section~\ref{sec:examples}, 
the true dimensions of slow manifolds are known and used. 
In applications when the size of the slow manifold dimension is not known, tools like diffusion maps~\cite{coifman2006diffusion} 
or proper orthogonal decomposition (POD)~\cite{berkooz1993proper} 
may be used to infer the size of the reduced dimension before 
applying the proposed FSNN.

\subsection{The Invertible Coupling Flow Network}
\label{sec:inn}

Now we will describe further details of the FSNN architecture.
The first component is the invertible function, $h(z;\epsilon)$, 
which represents the transformation into fast-slow coordinates.
The network is based on the invertible coupling flow network
developed in \cite{ardizzone2018}, 
where the authors closely follow the architecture 
provided by~\cite{dinh2016}. 
{This network is an example of an invertible neural network (INN).
In particular,  invertible coupling flow networks}
have been shown to be universal 
approximators to invertible functions in~\cite{teshima20}.
In this work, we modify the original network by {introducing} a 
regularization property which controls the Lipschitz constant of
both the forward and inverse map.
Such a regularization improves the robustness of the network
during training. 
The constant is a tunable parameter of the network that,
similar to the hidden dimensions of the network,
can be made arbitrarily large to increase the expressiveness of 
the network. Therefore, this modification is 
compatible with the universal approximation result 
derived in~\cite{teshima20}.
As such, we introduce the key ingredient of Lipschitz control.

{\definition 
The bi-Lipschitz affine transformation (bLAT) layer 
is given by the transformation, $f:\mathbb{R}^{N}\rightarrow\mathbb{R}^{M}$,
defined by
\begin{align}
f(x) = U \Sigma V^{\rm T} x + b,
\end{align}
where 
$U\in\mathbb{R}^{M\times r}$,
{$V\in\mathbb{R}^{N\times r}$},
$\Sigma\in\mathbb{R}^{r\times r}$,
$b\in\mathbb{R}^{M}$,
$r = \min(M, N)$,
$U$ and $V$ are orthogonal matrices satisfying 
$U^{\rm T} U = V^{\rm T} V = I$,
and 
\begin{align}
    \Sigma = {\rm diag}
    \left[\begin{array}{ccc} \sigma_1 & \dots & \sigma_r\end{array}\right],
\end{align}
where $\frac{1}{L} \le \sigma_i \le L$ for 
$i=1,\dots, r$ and some $L\ge 1$.
}

\medskip

\noindent {We remark that when $M=N$, the mapping is invertible.\label{def:blat}
\begin{align}
    f^{-1}(x) = V \Sigma^{-1} U^{\rm T} (x - b).
\end{align}}
Our new contribution is to include the
matrix term $U\Sigma V^T$, which is written as a 
singular value decomposition (SVD),
which exists for every matrix in $\mathbb{R}^{M\times N}$.
One advantage of the bLAT layer is the following theorem.

{\theorem The forward and 
inverse bLAT layers are both $L$-bi-Lipschitz.
}

\medskip

\noindent This follows from the fact that
$\frac{1}{L} \le \|\grad f(x)\|_2 \le L$ and 
$\frac{1}{L} \le \|\grad f^{-1}(x)\|_2 \le L$.
Another advantage of using the SVD is that,
by construction, the evaluation of 
$f(x)$ and $f^{-1}(x)$ have an equivalent computational cost.
In our implementation of the $M=N$ case,
the orthogonal matrices of the 
bLAT layers are parameterized using the matrix exponential 
of a skew-symmetric matrix input,
i.e. $U = {\rm expm}(S)$, where $S=-S^T$.
When $M\ne N$, the orthogonal matrices are parameterized
using the Householder factorization.

The second component of the invertible coupling flow network
is the affine coupling flow layer.
{We adapt the definition provided in Section 2.1
of~\cite{teshima20}
and introduce a similar regularization idea to the above.}

{\definition 
The affine coupling flow layer 
is defined by the transformation 
$y = {g}(x)$, $x, y \in\mathbb{R}^{M}$,
which is implicitly defined by 
\begin{align}
\left[\begin{array}{c}
              y_1 \\ y_2
            \end{array}\right]
            = 
\left[\begin{array}{cc}
              {\rm diag}(F(x_2)) & \\
                            & {\rm diag}(G(y_1))
            \end{array}\right]
                              \left[\begin{array}{c}
                                      x_1 \\ x_2
                                    \end{array}\right]
      +
      \left[\begin{array}{c}
              B(x_2) \\ C(y_1)
            \end{array}\right],
\end{align}
where 
$x = \left[\begin{array}{cc} x_1 & x_2 \end{array}\right]^{\rm T}$
$y = \left[\begin{array}{cc} y_1 & y_2 \end{array}\right]^{\rm T}$
$x_1, y_1\in\mathbb{R}^{d}$,
$x_2, y_2\in\mathbb{R}^{M-d}$,
$F, B: \mathbb{R}^{M-d}\rightarrow \mathbb{R}^{d}$,
$G, C: \mathbb{R}^{d}\rightarrow \mathbb{R}^{M-d}$, 
$1 \le d < M$
and the range of the output elements of $F$ and $G$ 
is constrained to 
$\left[\frac{1}{L},\, L\right]$ for some $L\ge 1$.
The inverse mapping, $x=f^{-1}(y)$, is implicitly defined by
\begin{align}
\left[\begin{array}{c}
              x_1 \\ x_2
            \end{array}\right]
            = 
\left[\begin{array}{cc}
              {\rm diag}(F(x_2))^{-1} & \\
                            & {\rm diag}(G(y_1))^{-1}
            \end{array}\right] \left(
                              \left[\begin{array}{c}
                                      y_1 \\ y_2
                                    \end{array}\right]
      -
      \left[\begin{array}{c}
              B(x_2) \\ C(y_1)
            \end{array}\right]\right).
\end{align}
}

\medskip

\noindent The functions $F, G, B, C$ 
are represented using feed-forward neural networks.
In order to constrain outputs to the range $\left[\frac{1}{L}, L\right]$, we apply the map $\theta_L: \mathbb{R}\rightarrow\left[\frac{1}{L}, L\right]$ defined by
$\theta_L(s) = \exp(\log(L) \tanh(s))$.
This constraint ensures invertibility of the
affine coupling layer.
The Lipschitz constant is a function of $L$ and
$\grad F, \grad G, \grad B, \grad C$ and can 
optionally be constrained using bLAT layers in the feedforward 
networks.

The invertible coupling flow network combines
bLAT layers and affine coupling flow layers into a network.

{\definition 
Let $f_1, \dots f_K:\mathbb{R}^{M}\rightarrow\mathbb{R}^M$ define a set of bLAT layers
and $g_1, \dots g_K:\mathbb{R}^{M}\rightarrow\mathbb{R}^M$ define a set of affine coupling flow
layers. The invertible coupling flow network,
$h:\mathbb{R}^{M}\rightarrow\mathbb{R}^{M}$
is defined
as the composition
\begin{align}
h = g_K \circ f_K \circ \dots \circ g_1 \circ f_1
\end{align}
and the inverse map is defined by
\begin{align}
h^{-1} = f_1^{-1} \circ g_1^{-1} \circ \dots \circ f_K^{-1} \circ g_K^{-1}.
\end{align}
}

\noindent The free parameter $L$ in each of the layers 
adjusts the complexity of the overall network and can
be used as an effective regularizer.
Despite the introduction of this regularization parameter, 
the network can be made arbitrarily expressive with 
the tuning of $L$, and therefore the universal approximation 
result from~\cite{teshima20} still holds true.

{\theorem The invertible coupling flow network is a universal 
approximator for invertible maps.
}

\medskip

\noindent In Theorem~\ref{thm:fsnn}, the invertible transformation 
is $\epsilon$-dependent. Therefore, we {incorporate} this dependence into
the network by including $\epsilon$ as a functional dependence in
$F, G, B, C$.

\subsection{The Negative Schur Form Network}
\label{sec:schur}


Next, a parameterization suitable for 
representing $T$ is discussed.
We present the negative Schur form network,
which is a differentiable parameterization
for Schur forms with eigenvalues 
to the left of the imaginary axis.
The Schur form network is composed of a 
feedforward neural network in which the outputs
are re-formatted into a negative Schur form matrix,
which is defined below.

{\definition The negative Schur form matrix
is given by the block matrix $T\in\mathbb{R}^{N_y\times N_y}$, where
\begin{align}
  T =
    \left[
    \begin{array}{cccc}
      T_{11} & T_{12} & \dots & T_{1p} \\
                     & T_{22} & \dots & T_{2p} \\
                     &                & \ddots & \vdots \\
                     &                &        & T_{pp}
    \end{array}\right]. 
  \label{eq:blockschur}
\end{align}
When $M$ is even, $p=M/2$ 
and $T_{ij}\in\mathbb{R}^{2\times 2}$ for $i=1,\dots p$, $j=i,\dots, p$.
When $M$ is odd, $p = (M+1)/2$,
$T_{ij}\in\mathbb{R}^{2\times 2}$ for $i=1,\dots p-1$, $j=i,\dots, p-1$,
$T_{ip}\in\mathbb{R}^{2}$ for $i=1,\dots, p-1$,
and $T_{pp}\in(-\infty, 0)$.
Each of 
the $\mathbb{R}^{2\times 2}$ blocks on the
diagonal are further parameterized using
the map
\begin{align}
\mathcal{M}(R, r, \theta, \phi) = 
  R \left[\begin{array}{cc}
            \cos(\theta) & \sin(\theta) \\
            -\sin(\theta) & \cos(\theta)
          \end{array}\right]
 +r \left[\begin{array}{cc}
            \cos(\phi) & \sin(\phi) \\
            \sin(\phi) & -\cos(\phi)
          \end{array}\right],
          \label{eq:2x2}
\end{align}
where the region of parameters is
restricted to 
\begin{align}
|R| \ge |r|, \quad
  R < 0, \quad
  \theta \in \left(-\frac{\pi}{2}, \frac{\pi}{2}\right),
  \quad
  {\phi \in \mathbb{R}}
  \label{eq:2x2restrict}
\end{align}
}

\medskip

\noindent Due to this parameterization,
we have the following theorem.

{\theorem The negative Schur form network
is a parameterization for the set of Schur 
forms with eigenvalues to the left of the
imaginary axis.
}

\medskip

\noindent This can be proven by 
showing that~\eqref{eq:2x2}--\eqref{eq:2x2restrict}
can smoothly represent both 
complex conjugate pairs and 
real pairs of eigenvalues with negative 
real part.
Equation~\eqref{eq:2x2} represents 
the $\mathbb{R}^{2\times 2}$
matrix as a 4-parameter family 
sum of a symmetric and
antisymmetric operator~\cite{Blinn1996ConsiderTL}.
Let $m = \frac{1}{2} {\rm Tr} (\mathcal{M}) = R \cos(\theta)$ and $p={\rm det}(\mathcal{M}) = R^2 - r^2$.
Then the eigenvalues of $\mathcal{M}$ 
are given by
\begin{align}
  \lambda_{\pm} = m \pm \sqrt{m^2 - p}.
\end{align}
The real parts of the eigenvalues are given by
\begin{align}
  {\rm Re}(\lambda_{\pm}) = m \pm \sqrt{\max(m^2 - p, 0)}.
\end{align}
Figure~\ref{fig:spectrumPlot} shows the regions of the $m{\rm-}p$ plane where both eigenvalues are less than zero
(green region), greater than zero (red region), complex (blue region), and real (outside of the blue region).
In order for ${\rm Re}(\lambda_{\pm}) < 0$, we must have
$m < 0$ and $p > 0$.
{In order to enforce these constraints, we enforce the restrictions
in~\eqref{eq:2x2restrict}.
We remark that for the first constraint, $R \cos(\theta) < 0$,
$R$ and $\cos(\theta)$ must have opposite signs.
If we focus on the region $\theta\in(-\pi/2, \pi/2)$
where $\cos(\theta)>0$, we would require $R<0$.
If the region of $\theta$ were instead shifted by $\pi$,
then we would require {$R>0$}. 
However, since $\mathcal{M}(R, r, \theta+\pi, \phi)
= \mathcal{M}(-R, r, \theta, \phi)$, these choices are equivalent.
}

{
\newcommand{\figWidth}{8cm}
\newcommand{\trimfig}[2]{\trimw{#1}{#2}{.0}{.0}{.0}{.0}}
\begin{figure}[thb]
\begin{center}
\begin{tikzpicture}[scale=1]
  \useasboundingbox (0,0) rectangle (13,5);  
  \draw( 2.6,-.4) node[anchor=south west,xshift=-15pt,yshift=-8pt] {\trimfig{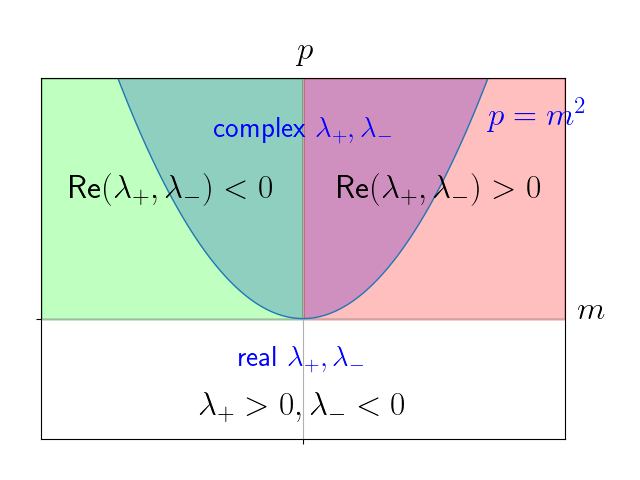}{\figWidth}};
\end{tikzpicture}
\end{center}
\caption{Plot showing the eigenvalues encountered in various regions of the $m$-$p$ plane, where
  $m = \frac{1}{2} {\rm Tr} (\mathcal{M}) = R \cos(\theta)$ and $p={\rm det}(\mathcal{M}) = R^2 - r^2$.
  The green region corresponds to linear stability, ${\rm Re}(\lambda_+, \lambda_-) < 0$, and the
  red region corresponds to linear instability, ${\rm Re}(\lambda_+, \lambda_-) > 0$.
  Above the parabola $p=m^2$ (in the blue region), eigenvalues are complex. Below the parabola (outside of the blue region),
  the eigenvalues are real.}
  \label{fig:spectrumPlot}
\end{figure}
}

In the implementation of the 
network,
each of the $\mathbb{R}^{2\times 2}$ blocks
have four tunable parameters,
$\bar{R},\, \bar{r},\, \bar{\theta},\, \bar{\phi} \in \mathbb{R}$,
that are 
mapped onto the restricted domain using
\begin{align}
  R = -|\bar{r}| e^{|\bar{R}|}, \quad
  r = \bar{r}, \quad
  \theta = \frac{\pi}{2} \tanh(\bar{\theta}), \quad
  \phi = \bar{\phi}.
\end{align}
The remaining off-diagonal terms 
are unrestricted in parameter space.
In order to achieve dependence with respect
to $x$, (i.e. $T(x)$), we introduce a 
feedforward network which takes $x$ as an
input and outputs
both the tunable parameters 
$\bar{R},\, \bar{r},\, \bar{\theta},\, \bar{\phi}$
for each diagonal block 
and the elements of the remaining blocks.
The resulting negative Schur form matrix is returned 
as the output of $T(x)$.

\subsection{The Low-Rank {bilinear map} Network}
\label{sec:bilinear}

Finally, we introduce an approximation 
for the {bilinear map} term, $B$.

{\definition The low-rank bilinear form
network $B:\mathbb{R}^{N_x}\times \mathbb{R}^{N_y}\rightarrow \mathbb{R}^{N_z}$ is defined by
\begin{align}
B(x, y) = 
\sum_{r=1}^{R} (C^{(r)} x) \odot (D^{(r)} y),
\label{eq:lowrank}
\end{align}
where 
$x\in\mathbb{R}^{N_x}$,
$y\in\mathbb{R}^{N_y}$,
$C^{(r)}\in\mathbb{R}^{N_z\times N_x}$,
$D^{(r)}\in\mathbb{R}^{N_z\times N_y}$,
$\odot$ represents the Hadamard product
(element-wise multiplication),
and $R$, $1\le R \le \min(N_x, N_y)$, represents 
the rank of approximation.
}

\medskip

\noindent Here, $R$ is a tunable parameter. Increasing $R$ 
gives more expressive power at the expense of 
computational cost. This network has the following 
property:

{\theorem The low-rank {bilinear map} network is a universal 
approximator to the set of bilinear forms.
}

\medskip

\noindent Consider the bilinear form, 
$\Bt:\mathbb{R}^{N_x}\times 
\mathbb{R}^{N_y}\rightarrow \mathbb{R}^{N_z}$.
{The tensor representation of $\Bt$ is written as} 
$\Bt_{ijk}$ where $\Bt_i\in\mathbb{R}^{N_x\times N_y}$
for $i=1,\dots, N_z$.
Let $p=\min(N_x, N_y)$, {then} we can use the singular value 
decomposition to represent the matrix $\Bt_i$ as
\begin{align}
    \Bt_i = \sum_{r=1}^{p} \sigma_i^{(r)} u_i^{(r)} (v_i^{(r)})^T,
\end{align}
where $\sigma_i^{(r)}$, $u^{(r)}$, $v^{(r)}$,
$r=1,\dots, p$, are the singular values, left singular 
vectors, and right singular vectors of $\Bt_i$, 
respectively, and 
$\sigma_1\ge\dots \ge \sigma_p\ge 0$.
The $i^{\rm th}$ term of $\Bt(x, y)$ 
is given by
\begin{align}
    x^T \Bt_i y &= 
    \sum_{r=1}^{p} \sigma_i^{(r)} 
    x^T u_i^{(r)} (v_i^{(r)})^T y, 
    = \sum_{r=1}^{p} \sigma_i^{(r)} \left(\sum_{j=1}^{N_x}  u_{ij}^{(r)} x_j\right)
    \left(\sum_{k=1}^{N_y} v_{ik}^{(r)} y_k\right),
\end{align}
We can find some $C\in\mathbb{R}^{N_z\times N_x}$
and $D\in\mathbb{R}^{N_z \times N_y}$ such that
\begin{align}
    x^T \Bt_i y &= 
    \sum_{r=1}^{p} \left(\sum_{j=1}^{N_x}  C_{ij}^{(r)} x_j\right)
    \left(\sum_{k=1}^{N_y} D_{ik}^{(r)} y_k\right).
\end{align}
This is equivalent to 
\begin{align}
    \Bt(x, y) &= 
    \sum_{r=1}^{p} (C^{(r)} x) \odot (D^{(r)} y).
\end{align}
The network uses a truncated version of this sum in {equation~\eqref{eq:lowrank}}
to approximate the action of the {bilinear map} $\Bt$.
The approximation error of the low-rank {bilinear map} network
is $\|B_i - \Bt_i\|_2 = \sigma_i^{(R+1)}$. This error
can be decreased by increasing $R$.

In our formulation, the operator $B$ has coefficients
which are functions of $x, y,$ and $\epsilon$.
We inject this functional dependence by making 
the elements of $C^{(r)}$ and $D^{(r)}$ outputs of a
feed-forward network with inputs $(x, y, \epsilon)$.

\subsection{Numerical Integration}

The Schur form, $T$ will have eigenvalues 
that can be of arbitrarily large magnitude.
Due to this, the dynamics of~\eqref{eq:fastslowthm} 
will be stiff.
Therefore, an Implicit-Explicit (IMEX) time integrator can be used to
solve~\eqref{eq:fastslowthm} efficiently.
The IMEX scheme has been applied to a neural ODE framework in~\cite{zhang2023semi}.
As shown below, an additional benefit of the proposed FSNN is to {eliminate} the need for iterative or direct linear solvers, thereby simplifying the training process.

Consider an IMEX Runge-Kutta method consisting of a diagonally implicit {method}
for treating the stiff linear term and an explicit method for treating the
nonlinear term. The method can be represented by the following Butcher tableau~\cite{pareschi}

\smallskip

\begin{center}
\begin{tabular}{cccc|cccc}
  $a_{11}$ &&&& \\
  $a_{21}$ & $a_{22}$ &&& $\alpha_{21}$\\
  $\vdots$ & $\vdots$ & $\ddots$ && $\vdots$ & $\ddots$\\
  $a_{s1}$ & $\dots$ & $\dots$ & $a_{ss}$ & $\alpha_{s1}$ & $\dots$ & $\alpha_{s,{s-1}}$ \\
  \hline
  $b_1$ & $\dots$ & $\dots$ & $b_s$ & $\beta_1$ & $\dots$ & $\beta_{s-1}$ & $\beta_{s}$
\end{tabular}
\end{center}

\medskip

\noindent Let $x_n \approx x(t_n)$, 
$y_n \approx y(t_n)$,
and $\dt = t_n-t_{n-1}$. 
Let $f(x, y, \epsilon)=B(x, y, \epsilon)(y, y) + \epsilon C(x, y, \epsilon) y$.
The intermediate stages are implicitly defined by
\bse
\begin{align}
    X^i &= x_{n-1} + \epsilon\dt \sum_{j=1}^{i-1} \alpha_{ij} {g(X^j, Y^j, \epsilon)}, \\
    Y^i &= y_{n-1} + \dt \sum_{j=1}^i a_{ij} {T(X^{j}) Y^j} + \dt \sum_{j=1}^{i-1} \alpha_{ij} {f (X^j, Y^j, \epsilon)} 
    \label{eq:tbs}
\end{align}
\ese
for $i=1, \dots, s$.
{In our notation, we take the above summation when $i=1$ to be 0 (i.e., $\sum_{j=1}^{0} (\cdot) = 0$). }
{Equation~\eqref{eq:tbs} can be rewritten to solve for $Y^i$}.
\begin{align}
     Y^i  &= (I - \dt a_{ii} T(X^{i}))^{-1} \left(y_{n-1} + \dt \sum_{j=1}^{i-1} \left(a_{ij} {T(X^{j}) Y^j} + \alpha_{ij} {f (X^j, Y^j, \epsilon)}\right)\right)
\end{align}
Note that the matrix given by $I - \dt  a_{ii} T(X^{i})$ is a
block 2-by-2 upper triangular matrix.
Due to the proposed architecture, this inverse can be handled efficiently using 
a backward substitution solve.
Since $T$ additionally has functional dependence, 
we treat this in a vectorization-friendly
matrix-free fashion, which avoids
storing multiple matrices $T(X^i)$ for each data pair 
in the batch set during training.
The solution at the next time step is obtained through
the update
\bse
\label{eq:additive}
\begin{align}
    x_n &= x_{n-1} + \epsilon\dt  \sum_{i=1}^s \beta_{i} {g(X^i, Y^i, \epsilon)}, \\
      y_n &= y_{n-1} + \dt \sum_{i=1}^s (b_i {T(X^{i}) Y^i} + \beta_{i} {f(X^i, Y^i, \epsilon)}).
\end{align}
\ese
In our implementation, we use the
two-stage second-order accurate IMEX strong-stability-preserving 
L-stable scheme defined in~\cite{pareschi}.

{During the prediction of the learned system, one can use the above IMEX time integrator to integrate the full system,
using small time steps to resolve the fast time scale. 
Alternatively, the proposed approach offers a unique capability to 
take larger time steps that are on the slow time scale
when only the slow dynamics is of interest.}
When the solution is initially on the slow manifold
($y=0$), the stages $Y^i$ for $i=1,\dots,s$ evaluate 
exactly to zero
and therefore the numerical scheme can be simplified 
to an integration of the slow variables given by
\bse
\label{eq:slowintegration}
\begin{align}
    x_n &= x_{n-1} + \epsilon\dt \sum_{i=1}^s \beta_{i} {g(X^i, 0, \epsilon)}, \\
    X^i &= x_{n-1} + \epsilon \dt \sum_{j=1}^{i-1} \alpha_{ij} {g(X^j, 0, \epsilon)}, \qquad i=1,\dots, s.
\end{align}
\ese
In this case, one can take larger time steps 
(e.g. $\Delta t \sim \frac{1}{\epsilon}$) and step over the fast dynamics formally. 
A hybrid scheme can also be used where 
the full system in~\eqref{eq:additive} can be integrated 
until $\|y\|$ is sufficiently small, after which 
slow integration in~\eqref{eq:slowintegration} is pursued
using larger time steps.
{We note that some of the first
works to use
Runge-Kutta based residual networks to learn 
equations from time series was pursued in~\cite{rico1992, rico1995nonlinear}.}

\medskip

\paragraph{Data-driven closure} 
A closure discovery procedure is implied through the above integration with large time steps. 
For a given full state vector $z_0$, 
the transformation $h$ first maps to the separate state of $(x_0, y_0)$. 
The state vector is then projected onto the closure of $y=0$; the system for $x$ is thus closed and evolved using~\eqref{eq:slowintegration}.
Finally, the evolved state of $x$ is lifted back to the full state $z$ using $h^{-1}$.


\section{Numerical Implementation}
\label{sec:numerical}

\begin{figure}
    \centering
    \includegraphics[width=1\textwidth]{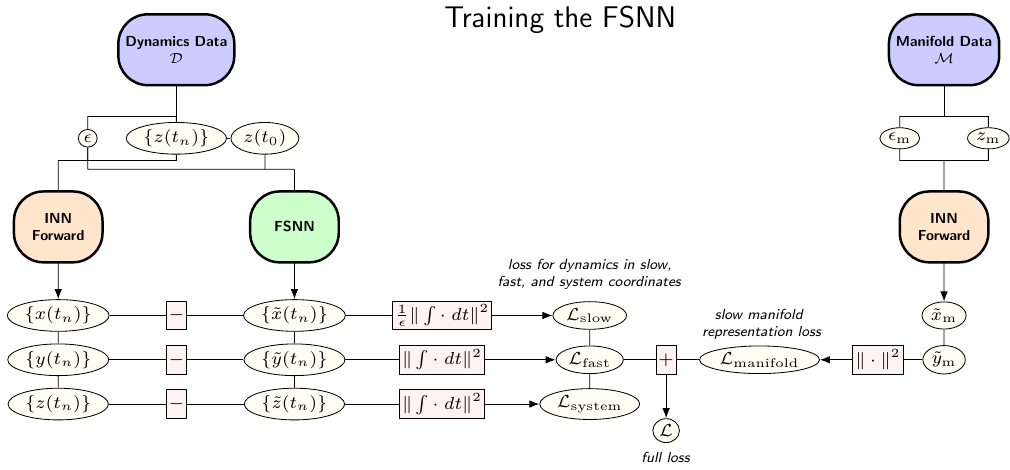}
    \caption{Schematic of training the fast-slow neural network (FSNN). The network is trained on dynamics data, $\mathcal{D}$,
    which contains examples of trajectories on both the fast and slow timescales
    and manifold data, $\mathcal{M}$, which contains examples 
    of solution vectors that lie on the slow manifold.
    {This data interacts both with the full network
    and the component involving the invertible neural network (INN).}
    The overall loss is a combination of the dynamics loss,
    which includes $L_2$ losses in the original coordinates
    ($\mathcal{L}_{\rm system}$), 
    fast coordinates ($\mathcal{L}_{\rm fast}$),
    and slow coordinates ($\mathcal{L}_{\rm slow}$),
    and the manifold loss ($\mathcal{L}_{\rm manifold}$),
    which is an $L_2$ loss of the fast coordinate, which 
    should map to $y=0$ on the slow manifold.
    }
    \label{fig:train}
\end{figure}

In this section, we describe the numerical implementation 
of the fast-slow neural network and its 
corresponding training algorithm.
The computation of the loss function and its gradient is summarized in Figure~\ref{fig:train}.
The training set consists of both trajectory data,
which describes the behavior of the dynamics, 
and manifold data, 
which describes the limiting behavior of the system.
The trajectory dataset, 
$\mathcal{D}$,
contains time series,
$\{(t_n^{(k)}, z_n^{(k)})\}_{n=1}^{N_t}$, 
and parameters $\epsilon^{(k)}$, 
for $k=1,\dots, N_d$, where $N_d$ is the number of 
trajectory data samples.
The manifold dataset, $\mathcal{M}$, 
contains multiple examples of the solution
on the slow manifold,
$z^{(k)}$, with corresponding parameters
$\epsilon^{(k)}$, for $k=1,\dots, N_m$.
Examples of manifold data can be obtained either 
experimentally by integrating a model for a sufficiently 
long time with a small $\epsilon$ or analytically 
if the equations permit an $\epsilon=0$ reduction
or asymptotic expansion. 
{For many examples using the zeroth order expansion with $\epsilon = 0$ is sufficient. Often obtaining a higher order expansion is nontrivial.}

The loss function is given by
\begin{align}
    \mathcal{L}
    = 
    \mathcal{L}_{\rm system}
    +\mathcal{L}_{\rm fast}
    +\mathcal{L}_{\rm slow}
    +\mathcal{L}_{\rm manifold}.
\end{align}
The first three terms are evaluated 
using the trajectory data set, $\mathcal{D}$,
and the last term is evaluated using the manifold
data set, $\mathcal{M}$.
For the dynamics-oriented loss terms, 
the fast-slow neural network is applied 
to initial conditions in $\mathcal{D}$
to obtain the model trajectories,
$\{(t_n^{(k)}, \tilde{x}_n^{(k)}, 
\tilde{y}_n^{(k)},
\tilde{z}_n^{(k)})\}$.
Additionally, the invertible coupling flow network 
is applied to the trajectory in $\mathcal{D}$
to augment the data
$\{(t_n^{(k)}, {x}_n^{(k)}, 
{y}_n^{(k)},
{z}_n^{(k)})\}$.
Let $L$ be an operator representing the 
approximation to the time integral of squared error, given by
\begin{align}
L(t, z, \tilde{z}) = 
    \sum_{n=0}^{N_t-1} \frac{1}{2}\left(
    (z_{n+1} - \tilde{z}_{n+1})^2
    + (z_{n} - \tilde{z}_{n})^2
    \right)(t_{n+1} - t_{n}).
\end{align}
Then we choose 
\bse
\begin{align}
    \mathcal{L}_{\rm system}
    &= \frac{1}{N_d}\sum_{k=1}^{N_d} L(t, z^{(k)}, \tilde{z}^{(k)}), \\
    \mathcal{L}_{\rm fast}
    &= \frac{1}{N_d}\sum_{k=1}^{N_d} L(t, y^{(k)}, \tilde{y}^{(k)}), \\
    \mathcal{L}_{\rm slow}
    &= {\frac{1}{ N_d}\sum_{k=1}^{N_d} \frac{L(t, x^{(k)}, \tilde{x}^{(k)})}{\epsilon^{(k)}}.}
\end{align}
\ese
Including the three terms individually 
encourages the network to learn the dynamics for
the original and fast-slow coordinates 
to an adequate level.

In practice, we split the trajectory dataset, $\mathcal{D}$,
into a fast-scale subset, with time steps chosen to 
smoothly capture the fast dynamics converging to the slow manifold,
and a slow-scale subset, with larger time steps 
chosen to capture the variation of the slow variables
after the dynamics sufficiently close to the slow manifold.
The trajectories in fast-scale subset are treated using 
the full integration in~\eqref{eq:additive}
while the trajectories from the slow-scale subset
are treated using the slow-manifold integration in~\eqref{eq:slowintegration}.

For the manifold loss term, the invertible coupling 
flow network is applied to examples in $\mathcal{M}$
to obtain fast-slow coordinates $(x^{(k)}, y^{(k)})$.
If the data is on the slow manifold, then 
$y^{(k)}$ should be zero. We therefore include
the following as a penalty term.
\begin{align}
    \mathcal{L}_{\rm manifold} = 
    \frac{1}{N_m} \sum_{k=1}^{N_m} (y^{(k)})^2
\end{align}

Our implementation of the fast-slow neural network uses
Jax~\cite{jax} as a backend. 
Auto-differentiation is used to compute gradients of the 
loss function with respect to the tunable parameters
and the Adam optimizer from the Optax~\cite{optax} package is 
used adjust the weights towards an locally optimal 
solution for the dataset.

\section{Examples}
\label{sec:examples}

In this section, we consider multiple example problems
to demonstrate the applicability of the FSNN to 
various problems, including multi-scale problems
arising in hydrodynamics and 
plasma physics.
{In training, we explored various combinations of meta parameters for each problem.
We pursued a policy to train several models in parallel and choose the model 
with the best performance. The meta parameters for the best models are reported in
each respective section.
We remark that though increasing the number of network parameters increases
the expressiveness of the model, it comes at a higher computational cost and increases the risk of overfitting.
}

\subsection{Visualization of the Attracting Slow Manifold Property}

A significant advantage of the FSNN is that, for any choice
of model parameters, the network will {represent a fast-slow dynamical system with}
an attracting
slow manifold that, 
due to Theorem~\ref{thm:fsnn}, will attract 
solutions for sufficiently small values of $\epsilon$.
The slow manifold is represented by the graph, $y=0$.
{Here we visualize the property of the FSNN without performing any training from dynamical system data.
Instead, we demonstrate the property by randomly generating the trainable weights  of 
FSNNs for a singularly perturbed dynamical system with $N_y=N_x=1$. The goal is to show that any initial condition for a randomly generated FSNN will naturally converge to its slow manifold. Note that for such a simple system, its slow manifold is known analytically. }

\begin{figure}[tbh]
    \centering
    \includegraphics[trim=4cm 0cm .5cm 0cm, width=1\textwidth]{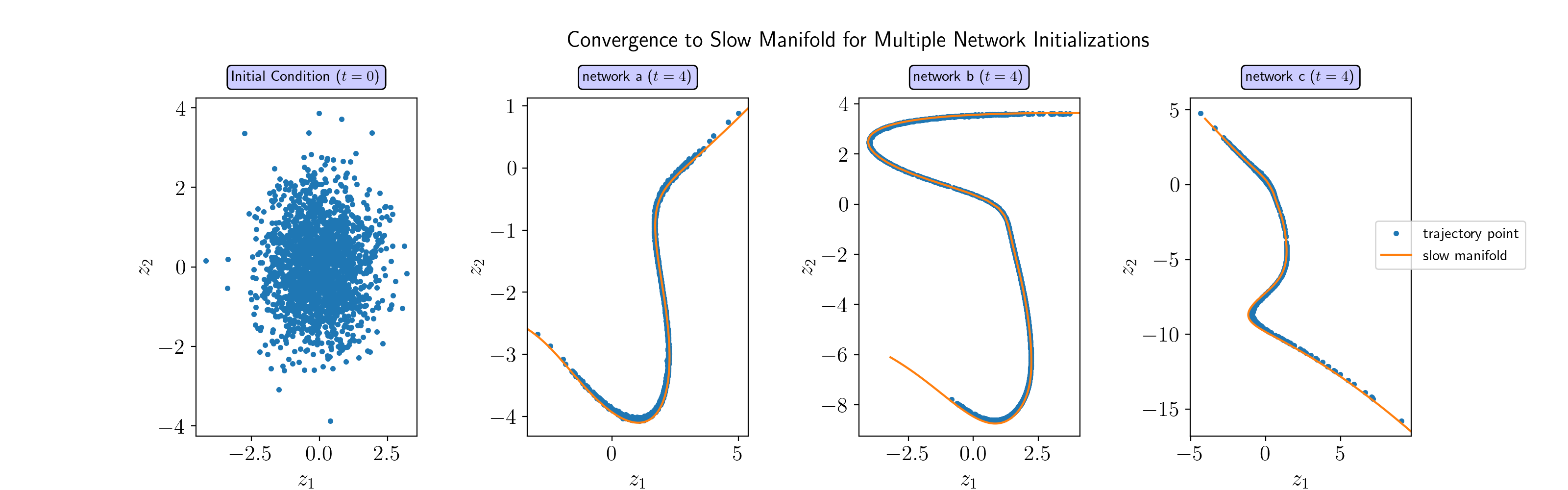}
    \caption{Demonstration of attracting slow manifold 
    property of the FSNN. For three different random initializations of the FSNN,
    initial conditions are drawn 
    randomly from a standard normal distribution
    and are evolved in time to a later time ($t=4$).
    The dynamics are shown to limit to the curve representing
    the slow manifold.}
    \label{fig:demo_manifold}
\end{figure}

{In Figure~\ref{fig:demo_manifold}}, network a, b, and c are representing 
different FSNNs with different weights and therefore
different slow manifolds.
Multiple initial conditions are 
drawn randomly from a standard normal distribution 
(left plot)
and are integrated using each FSNN to a later time $t=4$ {corresponding to $\epsilon=0$}. 
The trajectory points at $t=4$ (blue dots)
and the resulting slow manifold (orange curves) 
are plotted at this final time, demonstrating 
convergence of solutions to the slow manifolds {inherent to
each network}.

\subsection{Learning a Simple Attracting Manifold}

We demonstrate the ability of the FSNN to learn dynamics
for a simple test problem 
where the dynamics can be described by the following system of 
ODEs for $z_1, z_2 \in \mathbb{R}$.~
\bse
\label{eq:toy}
\begin{align}
    \frac{d z_1}{dt} &= -\epsilon (\sin(z_1) + z_2) \\
    \frac{d z_2}{dt} &= \lambda(z_1) (z_2 - \theta(z_1))
\end{align}
\ese
In this example, $z_1$ represents the slow variable and 
$z_2$ represents the fast variable. 
When $\epsilon=0$, the slow manifold is given by 
$\theta(z_1)$. 
In this example, we choose $\lambda$ and $\theta$ to be
\begin{align}
    \lambda(x) = -1 - \frac{1}{10} \cos(2 x),
    \qquad
    \theta(x) = 2 \tanh(x).
    \label{eq:evalsman}
\end{align}

The trajectory dataset, $\mathcal{D}$, 
consists of $10^3$ fast-scale trajectories, 
corresponding to initial conditions for $z_1, z_2$
randomly generated uniformly in $[-1, 1]$,
$\epsilon$ randomly generated log-uniformly 
in $[10^{-5}, 10^{-2}]$,
and integrated forward 5 timesteps using an adaptive 
RK45 solver with an error tolerance of $10^{-8}$.
The manifold dataset, $\mathcal{M}$, consists
of $10^3$ examples of data pairs $((z_1, z_2), \epsilon=0)$ 
on the graph $z_2 = \theta(z_1)$, where $z_1$
is randomly generated uniformly in $[-1, 1]$.

A FSNN is trained on this dataset using the Adam optimizer from
Optax. 
The INN consists of one outer layer, where the inner 
coupling flow layers consists of feedforward networks 
with one layer and a hidden dimension of 10.
Each of the feedforward neural networks on the RHS of the
neural ODE system consists of 10 layers and 
the {bilinear map} has a rank of 2.

\begin{figure}[tbh]
    \centering
    \includegraphics[trim=0cm 2cm 0cm 1cm, width=\textwidth]{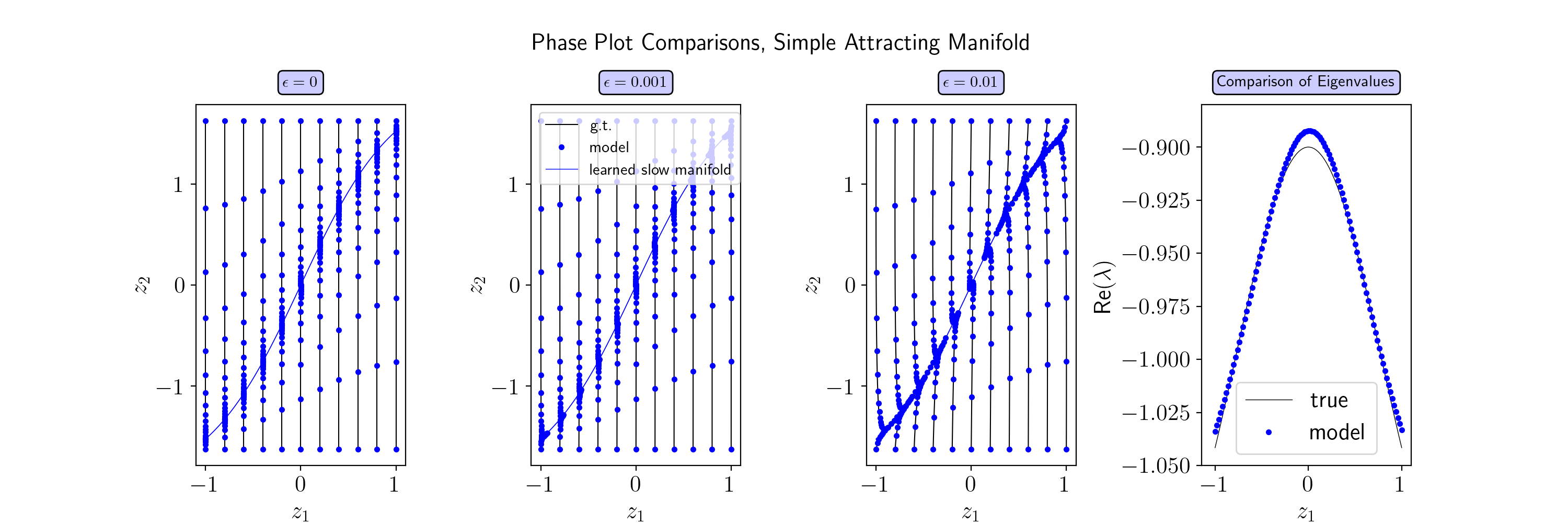}

\vspace{.5cm}
    
    \includegraphics[trim=0cm 0cm 0cm 0cm, width=\textwidth]{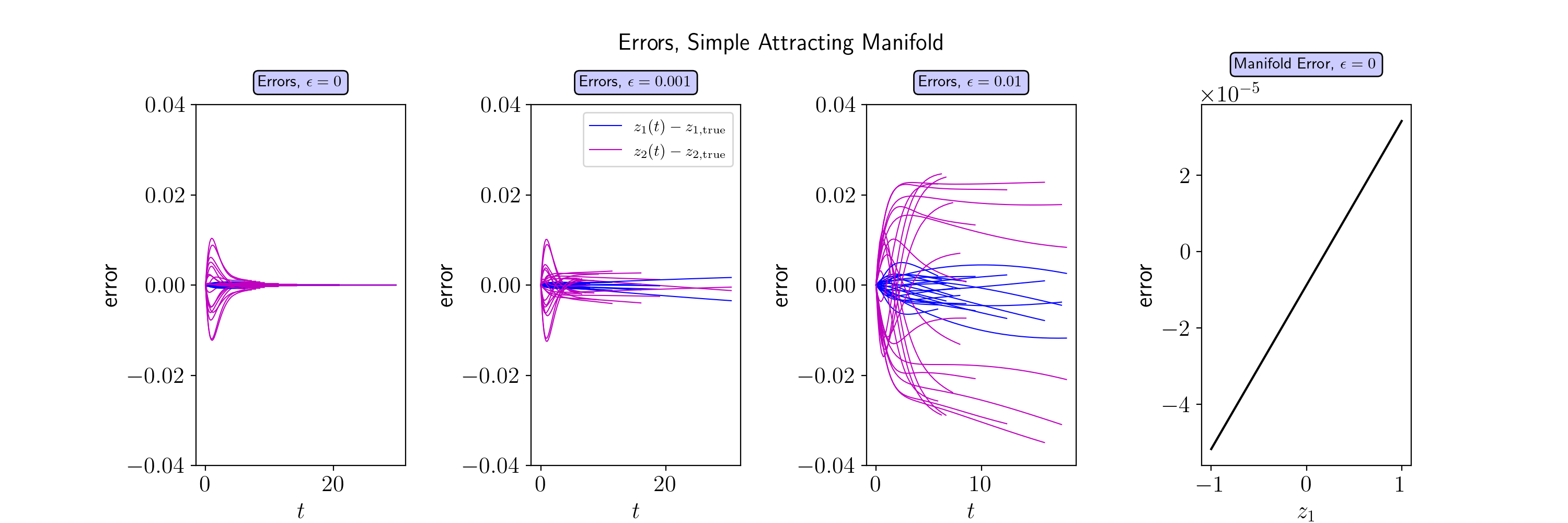}
    \caption{Results of learning the dynamics of the 
    simple attracting manifold 
    example described by~\eqref{eq:toy}.
    Left three plots on top row: phase portraits comparing the ground truth (g.t.) reference solution (bold curves) and that of the trained model (dots) for $\epsilon=0, 0.001, 0.01$.
    {Left three plots on bottom row: corresponding errors in time for each trajectory in the above phase portraits.}
    Rightmost plot on top row: comparison of the eigenvalue {of the negative Schur form as a function of $z_1$} learned by
    the FSNN with the ground truth eigenvalue {given by $\lambda$ in~\eqref{eq:evalsman}}.
    {Rightmost plot on bottom row: error between the learned slow manifold and the ground truth slow manifold given by $\theta$ in~\eqref{eq:evalsman}}.
    }
    \label{fig:simple_phase}
\end{figure}

Figure~\ref{fig:simple_phase} shows a comparison of 
the trained model and the reference solution in 
phase space for various choices of $\epsilon$.
{As can be seen by the corresponding error plots, the model shows reasonable agreement to the reference solution.}
In all of the phase plots, it is evident that the trajectories
eventually get trapped in a slow manifold.
In the rightmost plot {on the top row}, a comparison of the learned 
eigenvalues with that of the ground truth, $\lambda(z_1)$,
is shown. 
{Additionally in the rightmost plot on the bottom row, the error between the learned slow manifold and true slow manifold is shown.}
The results demonstrate that the FSNN is able to 
successfully recover details of the original operator.

\subsection{Linear Grad Moment System}
In this example, we show that the FSNN can be 
applied towards learning hydrodynamic slow manifolds.
We consider a linear Grad moment system~\cite{grad13, grad13kinetic}, defined by
\bse
\label{eq:gradmoment}
\begin{align}
\partial_t \rho &= - \epsilon \, \nabla \cdot \mathbf{u},\\
\partial_t \mathbf{u} &= - \epsilon \left(\nabla \rho + \nabla T + \nabla \cdot \sigma\right),\\
\partial_t T &= -\frac{2\epsilon}{3} \left(\nabla \cdot \mathbf{u} + \nabla \cdot \mathbf{q}\right),\\
\partial_t \mathbf{\sigma} &= -  \mathbf{\sigma} - \epsilon \left(2 \overline{\nabla \mathbf{u}} +\frac{4}{5} \overline{\nabla \mathbf{q}}\right) ,\\
\partial_t \mathbf{q} &= 
- \frac{2}{3} \mathbf{q}
-\epsilon \left(\frac{5}{2} \nabla T + \mathbf{\nabla} \cdot \mathbf{\sigma}\right) ,
\end{align}
\ese
where $\rho(\xv,t)$ is the density,
$\uv(\xv,t)$ is velocity,
$T(\xv,t)$ is temperature,
$\sigma(\xv,t)$ is the stress tensor,
$\mathbf{q}(\xv,t)$ is the heat flux,
and $\xv=(x, y, z)^T\in\mathbb{R}^3$ is the spatial coordinate.
The overline represents the symmetric, traceless part of 
a tensor, e.g.
\begin{align}
\overline{\av} &= \frac{1}{2} (\mathbf{a} + \mathbf{a}^T) - \frac{1}{3} \mathbf{I} \, \text{Tr} (\mathbf{a}).
\end{align}
We consider solutions that are 
one-dimensional ($\partial_y = \partial_z = 0$)
and $2\pi$-periodic in space, e.g.
\begin{align}
\rho(\xv, t) = \sum_{k=-\infty}^{\infty} \rho_k(t) e^{ikx}.
\end{align}
The Fourier coefficients are governed by the following system
of ODEs.
{\small
\begin{align}
    \frac{d }{d t}
    \left[\begin{array}{c}\rho_k \\ u_k \\ T_k \\ \sigma_k \\ q_k \end{array}\right]
    = 
    \left[\begin{array}{ccccc}
    0&-ik\epsilon&0&0&0\\
    -ik\epsilon&0&-ik\epsilon&-ik\epsilon&0\\
    0&-\frac{2}{3} ik\epsilon&0&0&-\frac{2}{3} ik\epsilon\\
    0&-\frac{4}{3} ik\epsilon&0&-1&-\frac{8}{15} ik\epsilon\\
    0&0&-\frac{5}{2}ik\epsilon&-ik\epsilon&-\frac{2}{3}
    \end{array}\right]
    \left[\begin{array}{c}\rho_k \\ u_k \\ T_k \\ \sigma_k \\ q_k \end{array}\right], \quad k=0, \pm 1, \dots
\end{align}
}

\remove{
Due to linearity, the equations for each Fourier coefficient
are independent.
The scalar version of the equations is:
\begin{align*}
\partial_t \rho_t &= -u_x,\\
\partial_t u &= -\rho_x -T_x - \sigma_x,\\
\partial_t T &= -\frac{2}{3} (u_x + q_x),\\
\partial_t \sigma &= -\frac{4}{3} u_x - \frac{8}{15} q_x -\frac{1}{\epsilon} \sigma,\\
\partial_t q &= -\frac{5}{2} T_x - \sigma_x - \frac{2}{3\epsilon} q.
\end{align*}
}

For this example, we truncate the Fourier series
to consider only the {$k=0, \pm 1$} models. 
Initial conditions are randomly generated for 
each Fourier coefficient
uniformly in the complex region $[-1, 1] \times [-i, i]$.
The solution is computed for 50 time steps 
both using 
$\Delta t_{\rm fast} = \frac{1}{4}$
and $\Delta t_{\rm slow} = \frac{1}{4 \epsilon}$.
{The model is trained using only the information 
from the first time step
and tested using the remaining of the trajectory.}
The slow-scale trajectories are offset 
$\Delta t_{\rm slow}$ from $t=0$ to ensure the solution
is sufficiently close to the slow manifold.
A FSNN is trained on this dataset using the 
Adam optimizer from Optax. 
The INN consists of one outer layer, where the inner 
coupling flow layers consists of feedforward networks 
with one layer and a hidden dimension of 40.
The feedforward neural networks on the RHS of the
neural ODE system each consist of 50 layers and 
the {bilinear map} has a rank of 1.
Considering each coefficient 
and its real and imaginary parts, the dimension of 
the learned model is $\mathbb{R}^{30}$, with
$N_x = 18$ and $N_y = 12$.

\begin{table}[tbh]
\caption{Comparison of ground truth eigenvalues of the Jacobian
of the fast equations with that
of the {negative Schur form} of 
the trained model for the Grad Moment 
System problem. \label{tab:grad}
}
\begin{center}
\begin{tabular}{|c|c|c||c|c|c|}
\hline
\multicolumn{6}{|c|}{Grad Eigenvalues} \\
\hline
 g.t. &   model &      error & 
 g.t. &   model &      error
 \\
 \hline
-1 & -0.9971 &   2.9e-03 &-$\frac{2}{3}$ & -0.6659 &   8.2e-04 \\
-1 & -0.9969 &   3.1e-03 &-$\frac{2}{3}$ & -0.6658 &   8.3e-04 \\
-1 & -0.9969 &   3.1e-03 &-$\frac{2}{3}$ & -0.6658 &   9.1e-04 \\
-1 & -0.9968 &   3.2e-03 &-$\frac{2}{3}$ & -0.6657 &   9.3e-04 \\
-1 & -0.9968 &   3.2e-03 &-$\frac{2}{3}$ & -0.6657 &   9.4e-04 \\
-1 & -0.9968 &   3.2e-03 &-$\frac{2}{3}$ & -0.6657 &   9.4e-04 \\
\hline
\end{tabular}
\end{center}
\end{table}

\begin{figure}[tbh]
    \centering
    \includegraphics[width=1\textwidth]{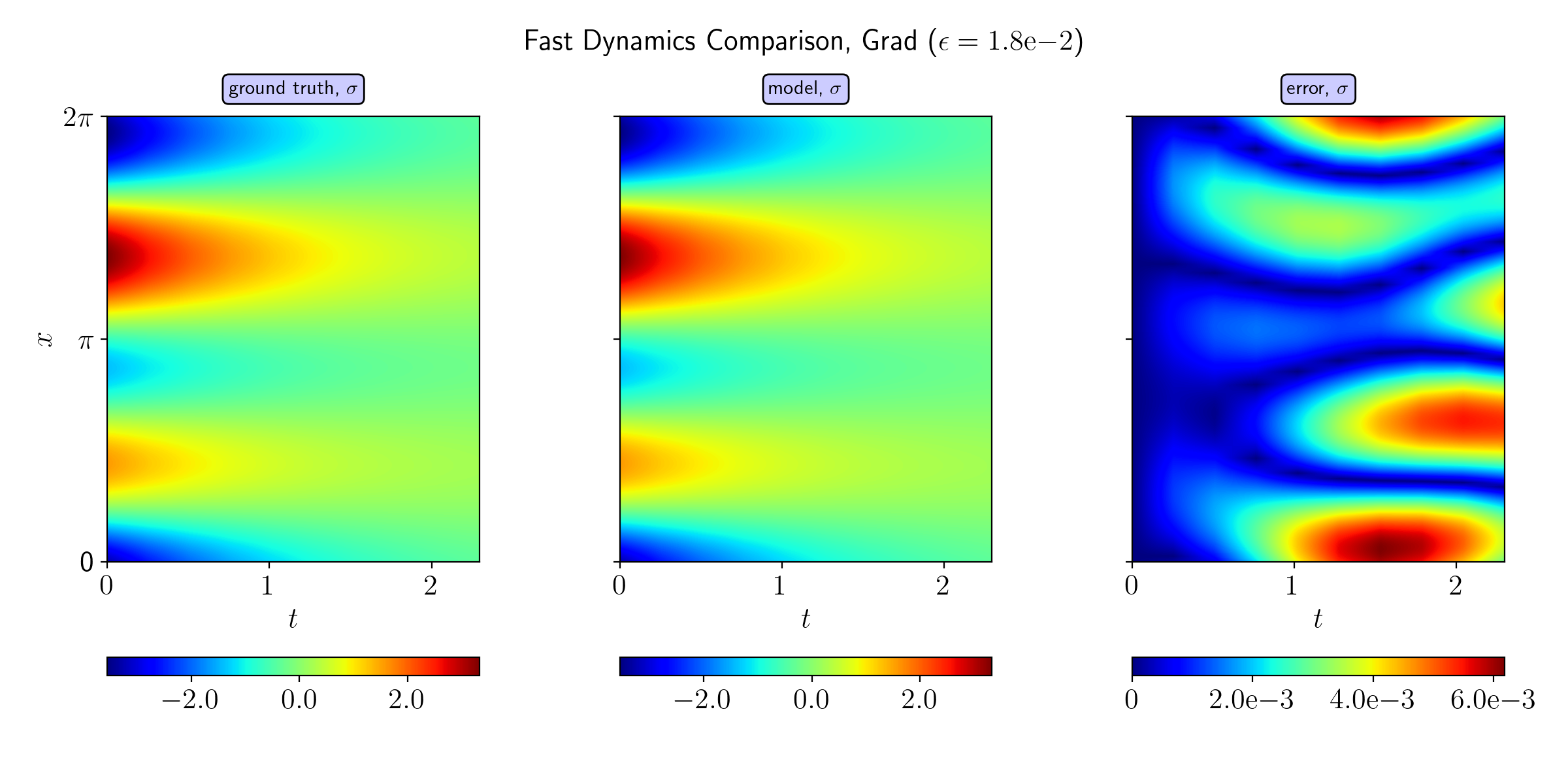}
    \includegraphics[width=1\textwidth]{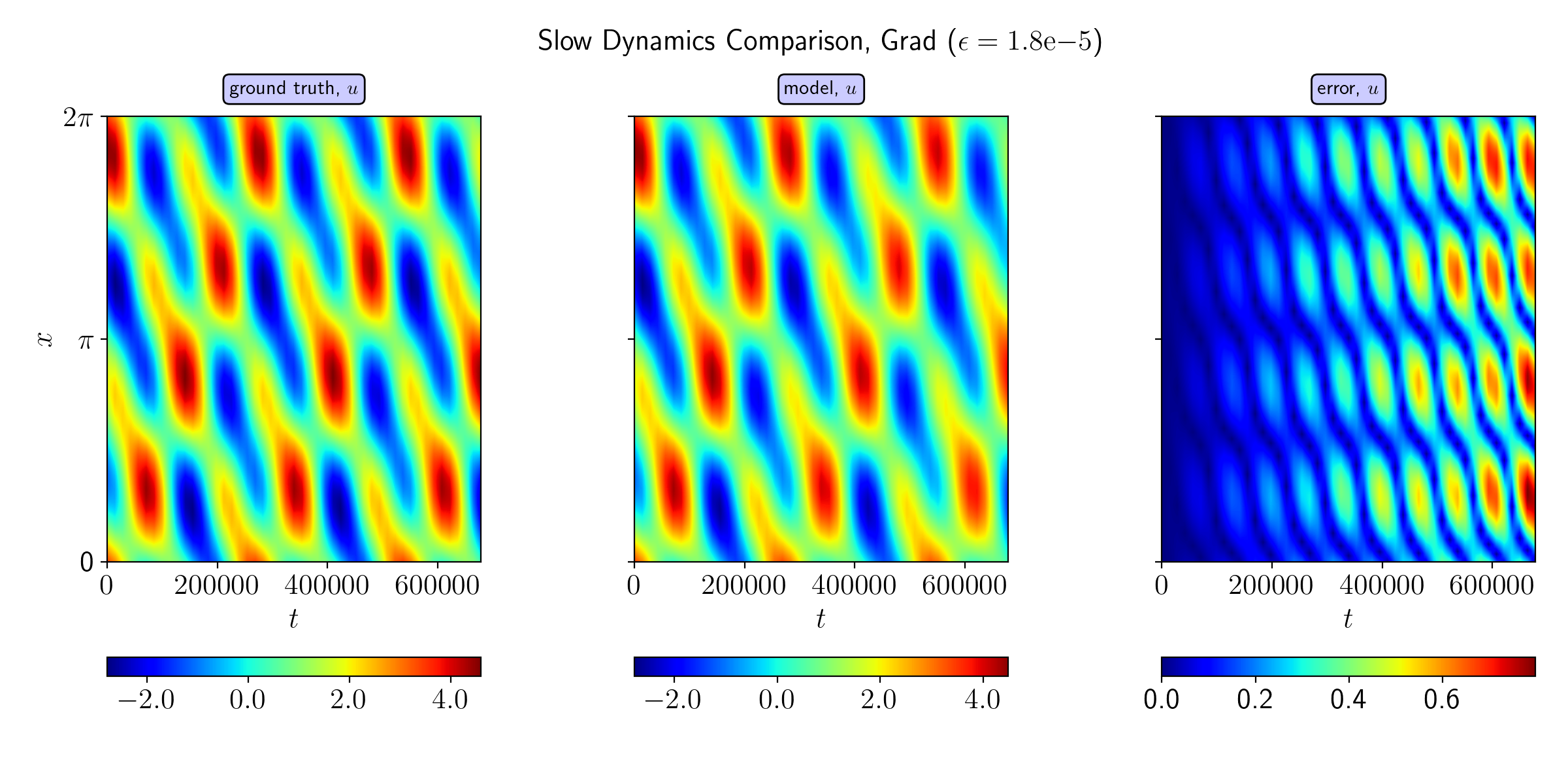}    
    \caption{Comparison of fast scale (top) 
    and slow scale (bottom) dynamics 
    between the ground truth 
    and trained model for the Grad moment system example.}
    \label{fig:gradcompare}
\end{figure}

Figure~\ref{fig:gradcompare} shows comparisons of 
the trained model and reference solutions on 
the fast time scale and slow time scale.
The results demonstrate that {reasonable accuracy can be achieved}
on both 
time scales.
{We observe that model for the slow dynamics 
is slightly off in phase which causes error to grow in time.}

Table~\ref{tab:grad} shows a comparison of the 
eigenvalues of the trained model with the analytical eigenvalues.
{Quantitative} agreement is observed, which implies that the model was successful in mapping the slow manifold approximately to the true slow manifold.

It is critical to stress that the slow dynamics comparison 
in~Figure~\ref{fig:gradcompare}
is done through integration on the slow manifold using~\eqref{eq:slowintegration}.
Specifically, the initial condition state vector $z_0$ is first transformed into the fast-slow coordinates $(x_0, y_0)$
and projected onto the learned slow manifold $(x_0, 0)$
before integrating the slow variable $x$
on the slow time scale only 
{and then transforming back into 
original coordinates using $h^{-1}$}.
This approach is significantly more 
efficient than integrating the full system
{using a standard numerical integrator.} 
The excellent accuracy of the predicted slow variable solution indicates a data-driven~\emph{closure} is learned through this procedure. 
The closure of~\eqref{eq:gradmoment} formally exists when $\epsilon$ approaches 0.
It is readily seen that the reduced system is a linearized Euler system.
{When $\epsilon=0$, the critical manifold is simply the center subspace.
When weak nonlinearity is added, there are center manifolds tangent to the center subspace.
In general these center manifolds are not unique, but there is a unique smoothest center manifold known as the spectral submanifold, see~\cite{haller2016nonlinear, cenedese}. 
In future work, it would be interesting to explore whether our 
methods can be adapted to ensure the learned slow 
manifold agrees with the spectral submanifold.}

\subsection{Multi-Scale Lorenz96 Equations}

In this section we consider the two-scale Lorenz96
equations, which was a model proposed in~\cite{lorenz96}
as a test bed model for multi-scale behavior in atmospheric 
dynamics.
{Multiscale techniques for this system have been explored
in many works, such as in~\cite{fatkullin04}.}
The original system is formulated as
\bse
\begin{align}
    \frac{d}{d\tau} X_k &= 
    -X_{k-1} (X_{k-2} - X_{k+1})
    - X_k
    + F
    - (hc/b) \sum_{j=1}^J Y_{j,k}, \\
    \frac{d}{d\tau} Y_{j,k} &=
    - cb Y_{j+1,k} (Y_{j+2,k} - Y_{j-1,k}) 
    - c Y_{j,k}
    + (hc/b) X_k,
\end{align}
\ese
for $j=1,\dots, J$ and $k=1,\dots, K$,
where the solutions satisfy periodicity 
($X_{k+K}=X_{k}$, $Y_{j,k+K}=Y_{j,k}$, $Y_{j+J,k}=Y_{j,k}$).
Let $c=b=1/\epsilon$, $h=1$, $Y=\epsilon y$, $\tau=\epsilon t$,
and $F=0$.
The system becomes 
\bse
\begin{align}
    \frac{d}{dt} x_k &= 
    -\epsilon x_{k-1} (x_{k-2} - x_{k+1})
    - \epsilon x_k
    - \epsilon^2 \sum_{j=1}^J y_{j,k}, \\
    \frac{d}{dt} y_{j,k} &=
    - y_{j+1,k} (y_{j+2,k} - y_{j-1,k}) 
    - y_{j,k}
    + x_k. \label{eq:lorenz96rhs}
\end{align}
\label{eq:lorenz96}
\ese
{Since the transformation is singular at
$\epsilon=0$, we study the equations in the
region $0 < \epsilon \ll 1$.}

Since the FSNN is limited to representing normally stable fast-slow
systems, we explore the region of applicability where~\eqref{eq:lorenz96}
is a normally fast-slow system.
The normally stable region of Lorenz96 is identified by the following theorem.
{\theorem The system of equations defined by $\eqref{eq:lorenz96}$ represent
a normally stable fast-slow system when $-\frac{1}{2}<x_k<\frac{8}{9}$, for $k=1,\dots, K$.
}

\smallskip

\noindent Using Definition~\ref{def:normfastslow}, we need to show that 
the Jacobian of the right-hand-side of~\eqref{eq:lorenz96rhs}
has eigenvalues strictly to the left of the imaginary axis.
Consider the following nonlinear equation $f(y;\xi)\in\mathbb{R}^{J}$
for $y\in\mathbb{R}^{J}$ and $\chi\in\mathbb{R}$.
\begin{align}
    f(y; \chi) = -y_{j+1} (y_{j+2} - y_{j-1}) - y_j + \chi.
    \label{eq:lorenznonlin}
\end{align}
The Jacobian of this equation, $D_y f$ is given by
\begin{align}
    (D_y f)_{jk} = \frac{\partial f_j}{\partial y_k}
    = -\delta_{jk} - \delta_{{j+1},k} (y_{j+2} - y_{j-1})
    - y_{j+1} (\delta_{j+2,k} - \delta_{j-1,k}),
\end{align}
where $\delta_{jk}$ represents the Kronecker delta.
We are interested in the eigenvalues of 
$D_y$ when $f(y; \chi)=0$.
Note that one solution to $f(y; \chi)=0$ is given by $y_j= \chi$ 
for $j=1,\dots,J$.
At this fixed point, the Jacobian becomes
\begin{align}
    (D_y f)_{jk}|_{y=\xi} 
    = -\delta_{jk} 
    - \chi (\delta_{j+2,k} - \delta_{j-1,k})
\end{align}
\remove{
\begin{align}
D_y f =
    \left[\begin{array}{cccccc}
    -1&0&-x&&&x \\
    x&-1&0&-x\\
    &\ddots&\ddots&\ddots&\ddots\\
    &&x&-1&0&-x\\
    -x&&&x&-1&0\\
    0&-x&&&x&-1
    \end{array}\right].
\end{align}}
This represents a circulant matrix, of which the eigenvalues
are given by
\begin{align}
    \lambda_j = -1 - \chi (e^{ij4\pi/J} - e^{-ij2\pi/J}),
    \qquad j=0,1,\dots, J-1.
\end{align}
Note that
\begin{align}
    {\rm Re}(\lambda_j) = -1 - \chi (\cos(4\xi_j) - \cos(2\xi_j)),
\end{align}
where $\xi_j= 4\pi j / J$.
We are interested in the region for $\chi$ where
${\rm Re}(\lambda_j)<0$ for all $j$. 
As $J\rightarrow\infty$,
we have the bound
$-\frac{9}{8} \le \cos(4\xi_j) - \cos(2\xi_j) \le 2$.
This implies that the region is given by
$-\frac{1}{2} < x < \frac{8}{9}$.
For a finite $J$, the bounds of $\cos(4\xi_j) - \cos(2\xi_j)$ are 
determined by the resolution of $\xi_j$, which is a discrete function of $j$.
This affords us to take a slightly wider region
of $x$ where
${\rm Re}(\lambda_j)<0$ for all $j$. 
For example, when $J=4$, 
$-\frac{1}{2} < x < 1$.

Homotopy continuation can be used on the polynomial system of equations 
defined by $f$ to reveal a high-order polynomial containing all
possible roots. Maple was used to evaluate all the roots 
for multiple values of $\chi$ and for $J=4, 5, 6$. 
Within the region $\chi\in(-\frac{1}{2}, \frac{8}{9})$,
the only real-valued solution is given by $y_j=\chi$ for $j=1,\dots,J$.
We conjecture that there is only a single real-valued solution 
to~\eqref{eq:lorenznonlin} when $\chi\in(-\frac{1}{2}, \frac{8}{9})$
for $J>6$.

We have shown that the multi-scale Lorenz96 equations
are normally-stable for a restricted range of $x_k$. 
However, for any $x_k$, $y_{j,l}$, $\epsilon>0$,
the system is energy stable, which implies that,
for any initial condition,
the energy will decay in time
and the dynamics will eventually be able to be described by 
a normally fast-slow system for some time $t>0$.

{\theorem The equations defined by~\eqref{eq:lorenz96}
are energy-stable for $\epsilon>0$.
}

\smallskip

\noindent Consider again the original system. We multiply the 
equations by $x_k$ and $y_{j, k}$, respectively.
Let $(\cdot)' = \frac{d}{dt} (\cdot)$, then
\bse
\begin{align}
    x_k x_k' &= 
    - \epsilon x_k x_{k-1} (x_{k-2} - x_{k+1})
    - \epsilon x_k^2 
    - \epsilon^2 \sum_{j=1}^{J} x_k y_{j,k} \\
    y_{j,k} y_{j,k}' &= 
    - y_{j,k} y_{j+1,k} (y_{j+2,k} - y_{j-1,k})
    - y_{j,k}^2 + x_k y_{j,k}.
\end{align}
\ese
We collect the terms and take the sum over the discrete variables.
\begin{align*}
    \frac{1}{2}\sum_{k=1}^K (x_k^2)' &= 
    - \epsilon \sum_{k=1}^K x_k x_{k-1} x_{k-2} 
    + \epsilon \sum_{k=1}^K x_{k+1} x_k x_{k-1} 
    - \epsilon \sum_{k=1}^K x_k^2 
    - \epsilon^2 \sum_{k=1}^K \sum_{j=1}^{J} x_k y_{j,k} \\
    \frac{1}{2} \sum_{j=1}^J (y_{j,k}^2)' &= 
    - \sum_{j=1}^J y_{j,k} y_{j+1,k} y_{j+2,k}
    + \sum_{j=1}^J y_{j-1,k} y_{j,k} y_{j+1,k}
    - \sum_{j=1}^J y_{j,k}^2 
    + \sum_{j=1}^J x_k y_{j,k}
\end{align*}
Due to periodicity, the cubic terms evaluate to zero.
The remaining terms give us the following energy estimate.
\begin{align}
    \frac{1}{2}\sum_{k=1}^K (x_k^2)'
    + \frac{1}{2} \epsilon^2 \sum_{k=1}^K \sum_{j=1}^J (y_{j,k}^2)'
    = - \epsilon \sum_{k=1}^K x_k^2 
    - \epsilon^2 \sum_{k=1}^K \sum_{j=1}^{J} y_{j,k}^2  < 0
\end{align}

We train a FSNN to learn the multi-scale Lorenz96 equations
from data {choosing $J=K=4$}.
The trajectory dataset, $\mathcal{D}$, 
consists of $10^4$ fast-scale trajectories,
corresponding to initial conditions for $x_k$
randomly generated uniformly in $[-\frac{1}{2}, 1]$
and $y_{j,k}$ randomly generated from a normal distribution,
$\epsilon$ randomly generated log-uniformly 
in $[10^{-5}, 10^{-2}]$,
and integrated forward 10 timesteps using an adaptive 
RK45 solver with error tolerance $10^{-8}$.
The manifold dataset, $\mathcal{M}$ consists
of $10^4$ examples of data pairs $((x_k, y_{j,k}), \epsilon=0)$ 
on the graph $y_{j,k} = x_k$, where $x_k$
is randomly generated uniformly in $[-\frac{1}{2}, 1]$.

A FSNN is trained on this dataset using the Adam optimizer from
Optax. 
The INN consists of one outer layer, where the inner 
coupling flow layers consists of feedforward networks 
with one layer and a hidden dimension of 40.
The feedforward neural networks on the RHS of the
neural ODE system each consist of 40 layers and 
the {bilinear map} has a rank of 4.

\begin{figure}[tbh]
    \centering
    \includegraphics[width=.9\textwidth]{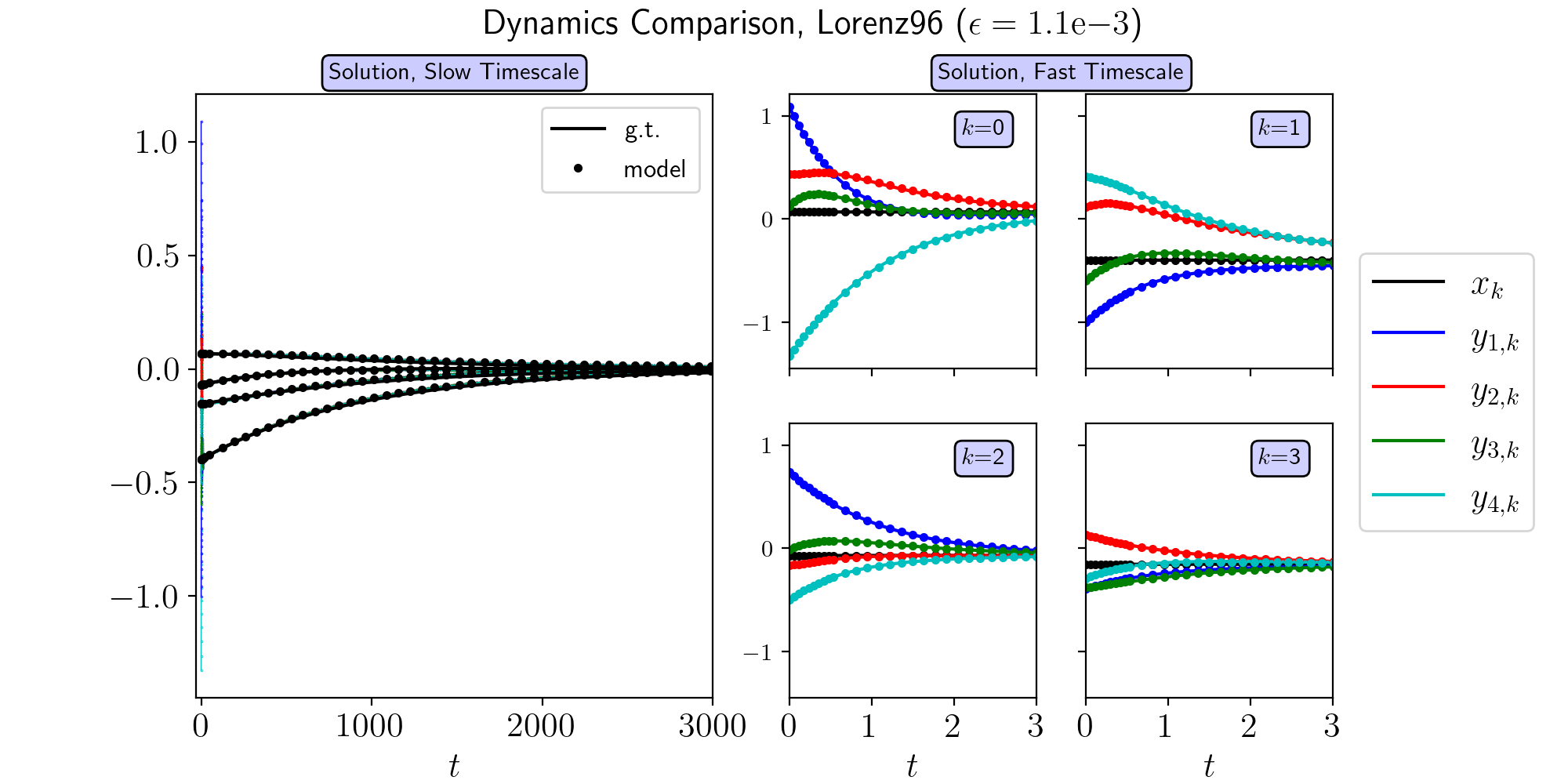}
    \includegraphics[width=.9\textwidth]{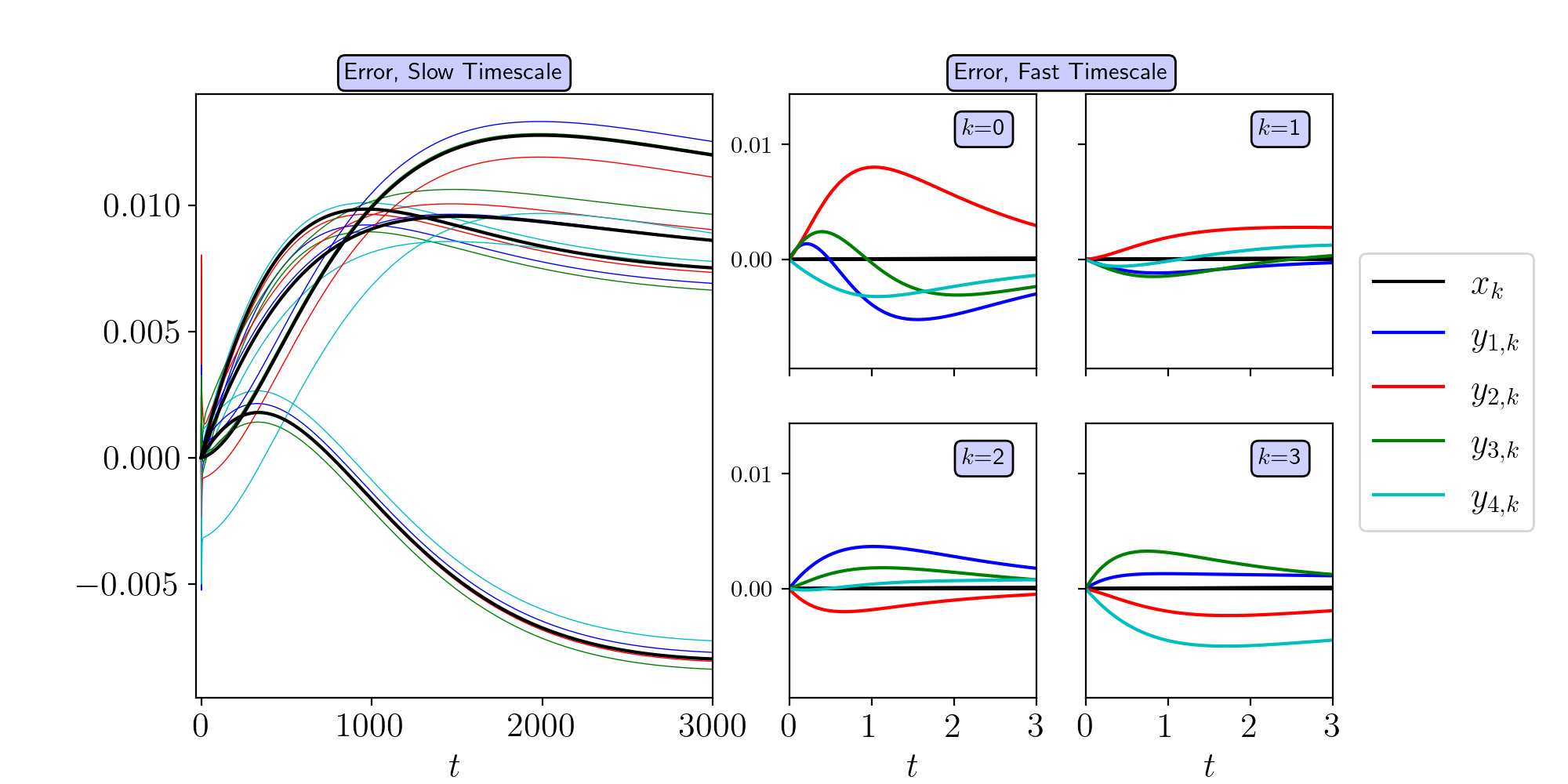}
    \caption{The top five plots are trajectory plots comparing the reference solution computed using
    RK45 (bold curves) and that of the trained model (dotted markers) for an initial condition with $\epsilon=1.1 \cdot 10^{-3}$ for the 
    multi-scale Lorenz96 example described by the dynamics 
    in~\eqref{eq:lorenz96}. 
    The bottom five plots show the corresponding errors
    as a function of time.
    The plots in the left column show the solution and 
    error on the slow time scale ($t\in [0, 3000]$)
    while the plots in the right column 
    zoom in on the fast time scale 
    ($t\in [0, 3]$).
    }
    \label{fig:lorenz_plots}
\end{figure}

\begin{figure}[tbh]
    \centering
    \includegraphics[width=.5\textwidth]{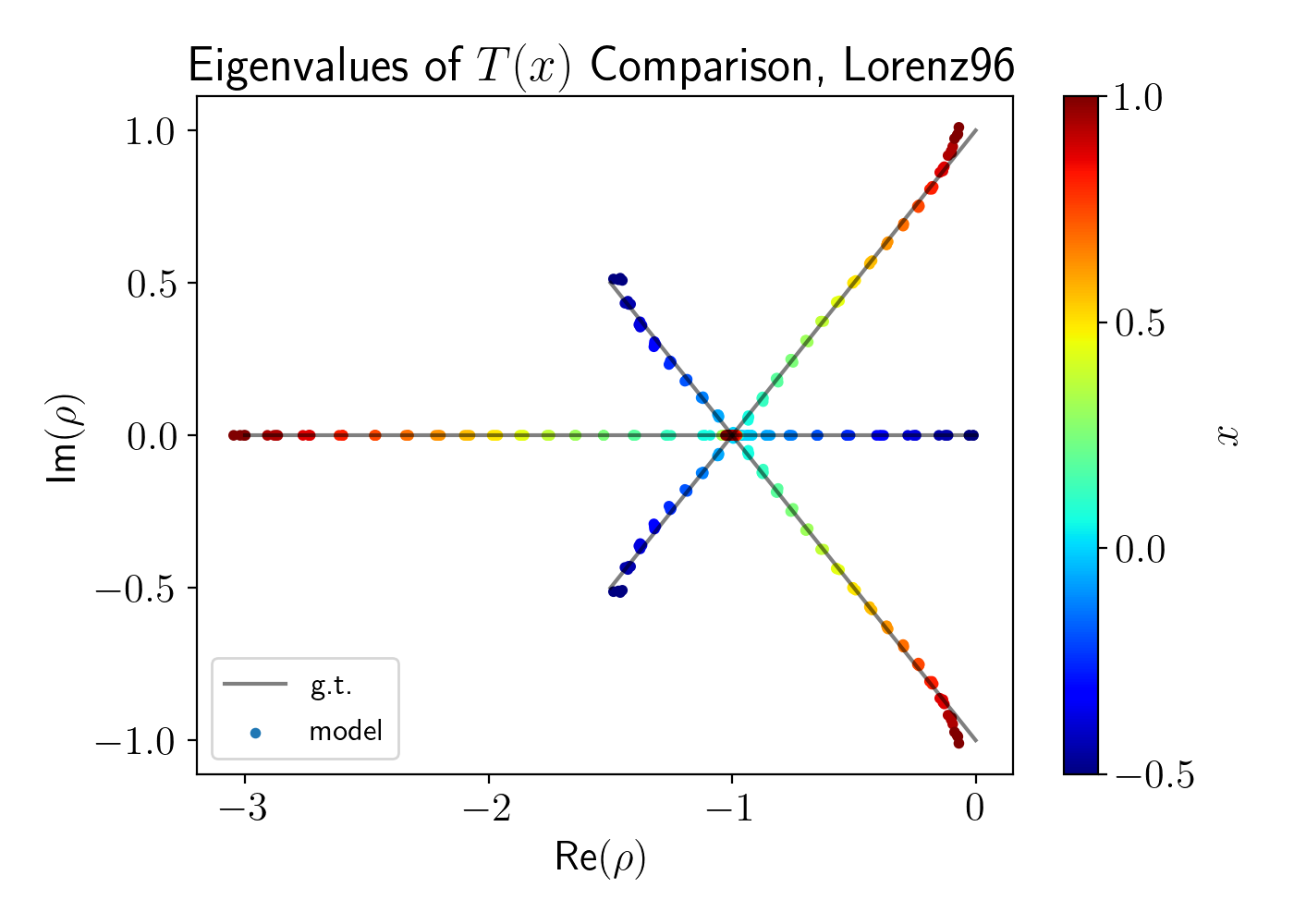}
    \caption{Comparison of the {ground truth eigenvalues of the Jacobian of the fast dynamics} (bold lines) with {that of the negative Schur form}
    learned from the trained model (dots)
    for the 
    multi-scale Lorenz96 example described by the dynamics 
    in~\eqref{eq:lorenz96}}
    \label{fig:lorenz_eig}
\end{figure}

Figure~\ref{fig:lorenz_plots} shows a comparison of 
the trained model and the reference solution for 
an example initial condition with $\epsilon=1.1 \cdot 10^{-3}$.
The model shows quantitative agreement to the reference solution.
{Although the globally attracting solution is the zero solution, the error remains away from 0. 
This can be explained 
by approximation error of the model
and can be improved through including more data in the
training set.
}
Figure~\ref{fig:lorenz_eig} shows a comparison of the 
eigenvalues of the trained model with the analytical eigenvalues.
Agreement is observed, implying that the model was successful 
in mapping the slow manifold approximately to the true slow manifold.
In this example, the time interval of all the training data is much shorter than $t=3$.
Nevertheless, the prediction can be done up to $T=3000$ or longer
with excellent accuracy.

\subsection{Abraham-Lorentz Equations}

\newcommand{\nooverline}[1]{#1}

In this section, we consider the Abraham-Lorentz equations~\cite{burby2020}, 
which model the radiation reaction in the case of a single electron with 
charge $e$ and mass $m$ moving through a static external magnetic field $\mathbf{B}(\mathbf{x})\in\mathbb{R}^{3}$. The equations are given by
\bse
\begin{align}
    \frac{2}{3} \frac{e^2}{c^3} \frac{d\mathbf{a}}{d\tau} &= m\mathbf{a} - \frac{e}{c} 
    \mathbf{v}
    \times \mathbf{B}(\mathbf{x}),\\
    \frac{d\mathbf{v}}{d\tau} &= \mathbf{a},\\
    \frac{d\mathbf{x}}{d\tau} &= \mathbf{v},
\end{align}
\ese
where $\mathbf{x}(\tau), \mathbf{v}(\tau), \mathbf{a}(\tau)\in\mathbb{R}^3$
are the position, velocity, and acceleration, respectively,
and $c$ is the speed of light.
We rescale these equations following the approach in~\cite{burby2020}. Time is scaled by the observer time scale $T$ as $\tau = T \overline{\tau}$, space by the observer length scale $\mathbf{x} = L \overline{\mathbf{x}}$, velocity as $\mathbf{v}=(L/T)\overline{\mathbf{v}}$, the magnetic field as $\mathbf{B}(\mathbf{x}) = B_0 \overline{\mathbf{B}}(\mathbf{x})$, and acceleration as $a = (L/T) (|e|B_0)(mc)^{-1} \overline{\mathbf{a}}.$
We reformulate the rescaled Abraham-Lorentz equations 
in terms of the nondimensional variables and drop the overline
for notational convenience.
\bse
\begin{align}
\left(\frac{r_0}{cT} \right)\frac{2}{3} \frac{d \nooverline{\mathbf{a}}}{d\nooverline{\tau}} & = \mathbf{a} - \zeta \mathbf{v} \times \mathbf{B}(\mathbf{x}),\\
\frac{d \mathbf{v}}{d\nooverline{\tau}} &= (\omega_c T) \mathbf{a},\\
\frac{d \mathbf{x}}{d\nooverline{\tau}} & = \mathbf{v},
\end{align}
\ese
where we have introduced the cyclotron frequency $\omega_c = |e|B_0/(mc)$, the classical electron radius $r_0 = e^2/(m c^2)$, and the sign of the charge $\zeta \in \{-1, +1\}$,
where $\zeta=-1$ corresponds to electrons
and $\zeta=+1$ to positrons.
The dimensionless parameters 
\begin{align}
\gamma_R = \frac{r_0}{cT}, \qquad 
\gamma_B = \frac{1}{\omega_c T},
\end{align}
represent the ratio of the electron size to the distance light travels during a time interval $T$, and ratio of the cyclotron period to $T$, respectively. 
The system can be transformed into a two-scale fast-slow system
with the choice $\gamma_B = 1$ and 
$\gamma_R~=~\epsilon~\ll~1$. 
We scale time through the transformation 
$\tau = \frac{2}{3} \epsilon \nooverline{t}$
to reveal a normally-unstable fast-slow system given by
\bse
\label{eq:forward}
\begin{align}
\frac{d \mathbf{a}}{dt} & = \mathbf{a} - \zeta \mathbf{v} \times \mathbf{B}(\mathbf{x}),\\
\frac{d \mathbf{v}}{dt} &=  \epsilon \frac{2}{3} \mathbf{a},\\
\frac{d \mathbf{x}}{dt} & = \epsilon \frac{2}{3} \mathbf{v},
\end{align}
\ese
We will refer to~\eqref{eq:forward} as the 
\emph{forward-time Abraham-Lorentz system}.
A normally-stable fast-slow system can be 
obtained by reversing time using $t_- = -t$,
\bse
\label{eq:reverse}
\begin{align}
\frac{d \mathbf{a}}{dt_-} & = -\mathbf{a} + \zeta \mathbf{v} \times \mathbf{B}(\mathbf{x}),\\
\frac{d \mathbf{v}}{dt_-} &=  -\epsilon \frac{2}{3} \mathbf{a},\\
\frac{d \mathbf{x}}{dt_-} & = -\epsilon \frac{2}{3} \mathbf{v}.
\end{align}
\ese
We refer to~\eqref{eq:reverse} as the 
\emph{reverse-time Abraham-Lorentz system}.
In both systems, $\mathbf{a}$ is the fast variable 
and $\mathbf{v}, \mathbf{x}$ are the slow variables.
In reverse time, the acceleration, $\mathbf{a}$
will be attracted to the slow manifold, given by
\begin{align}
\mathbf{a}^*(\mathbf{x}, \mathbf{v}) = 
\zeta \mathbf{v} \times \mathbf{B}(\mathbf{x}) + 
\epsilon \frac{3}{2} \mathbf{v} \times \mathbf{B}(\mathbf{x}) \times \mathbf{B}(\mathbf{x}) + \mathcal{O}(\epsilon^2).
\label{eq:expansion}
\end{align}
Though the reverse-time system is stable on the fast time scale,
the system exhibits an instability tangential to the slow manifold
for which solutions blow up on the slow time scale.
This is a more challenging test than the previous ones.

It is of interest to obtain forward-time trajectories 
of electrons along the slow manifold. 
A standard numerical integration of~\eqref{eq:forward}
will result in unstable trajectories on the fast
time scale since any numerical error will push trajectories
off of the slow-manifold, resulting in exponential growth
of the acceleration.
We instead apply FSNNs to learn the Fenichel normal form 
of the system in \emph{reverse time}. 
We then obtain forward-time trajectories 
along the learned slow manifold by performing 
an integration of the slow variables along $y=0$ using~\eqref{eq:slowintegration}.
This test is therefore a more practical application of a data-driven closure.

\remove{
{Here we reversed the time direction as well as scaling it. The reason for changing our time direction is that the Abraham-Lorentz equations are stable in backwards time.}
Here we assume $\epsilon \ll 1$. The formal slow manifold is given by 
{check the numbers on as they change frequently}
}

We train {an} FSNN to learn the Abraham-Lorentz equations
from data.
The trajectory dataset, $\mathcal{D}$, consists of $10^4$ trajectories, corresponding to initial conditions randomly generated in the region, $\mathcal{R}$, specified by
\begin{align*}
\mathcal{R} = \left\{\mathbf{x}, \mathbf{v}, \mathbf{a}: \quad
2 \le \nooverline{r} \le 4,\quad 
|z| \le 1, \quad
\|\mathbf{v}\| \le 1,\quad
\|\mathbf{a} - \mathbf{a}^*(\mathbf{x}, \mathbf{v})\| \le 1\right\},
\end{align*}
where $\nooverline{r}=\sqrt{\nooverline{x}_1^2 + \nooverline{x}_2^2}$ and $\nooverline{z} = \nooverline{x}_3$.
$\epsilon$ is randomly generated log-uniformly 
in the range $[10^{-5}, 10^{-2}]$.
The trajectories are 
integrated in reverse time using equations~\eqref{eq:reverse}
using an adaptive 
RK45 solver set to an error tolerance of $10^{-8}$.
The first 40 time steps batched in sequences of 
length 10 are used as fast-scale training data.
These trajectories are also interpolated
over a slower time scale with 
$\Delta t= 0.025 / \epsilon$ 
and $10^4$ sequences of length 10 are chosen at random
for slow-scale training data.
The manifold dataset, $\mathcal{M}$ consists
of $10^4$ examples of data pairs $((\mathbf{x}, \mathbf{v}, \mathbf{a}), \epsilon=0)$ 
on the graph $\mathbf{a} = \mathbf{v} \times \mathbf{B}(\mathbf{x})$, where $\mathbf{x}, \mathbf{v}$
are randomly generated in the region $\mathcal{R}$.

The magnetic field is chosen to be 
$\mathbf{B}  = B_r \mathbf{e}_r + B_\phi \mathbf{e}_\phi + B_z \mathbf{e}_z$,
where $\mathbf{e}_r, \mathbf{e}_\phi, \mathbf{e}_z$ represent basis vectors for 
cylindrical coordinates, and
\begin{align}
B_r & = -\frac{z}{q(r,z) r}B_0, \qquad
B_{\phi} = \frac{R_0}{r} B_0,\qquad
B_z = \frac{(r-R_0)}{q(r,z) r} B_0,
\end{align}
where $q(r,z) = q_0 + q_2((r-R_0)(r-R_0) + z^2)/a^2$ and $q_0 = 1.2$, $q_2 = 2.8$, $a = 1.5$, $R_0 = 3.0$, $Z_0 = 3.0$, and $B_0 = 1.0$.
This field is an idealized equilibrium of a tokamak, where the field is physical inside a torus with major radius of $R_0$ and minor radius of $a$.
Here $q$ is called a safety factor
\cite{liu2021parallel} and 
the poloidal magnetic flux function is $\psi = \frac{B_0 a^2}{2 q_2} \ln(q).$
Under this definition, the given field can be rewritten as
\begin{align}
    \Bv =\frac{1}{R} \nabla\psi\times\mathbf{e}_\phi + \frac{F_0}{R} \mathbf{e}_\phi.
\end{align}
where the poloidal field is a constant of $F_0=R_0 B_0$.
This test thus models a charged particle moving through a tokamak with radiation reaction.

A FSNN is trained on this dataset using the Adam optimizer from
Optax. 
The INN consists of one outer layer, where the inner 
coupling flow layers consists of feedforward networks 
with one layer and a hidden dimension of 50.
The function $g$ on the right-hand-side of the neural ODE
uses a single-layer feedforward neural network with 
hidden dimension 100 and is added to a bilinear form
network with rank 1.
The remaining 
feedforward neural networks on the right-hand-side of the
neural ODE system each consist of 50 layers and 
the {bilinear map} has a rank of 4.

\begin{figure}[tbh]
    \centering
    \includegraphics[width=\textwidth]{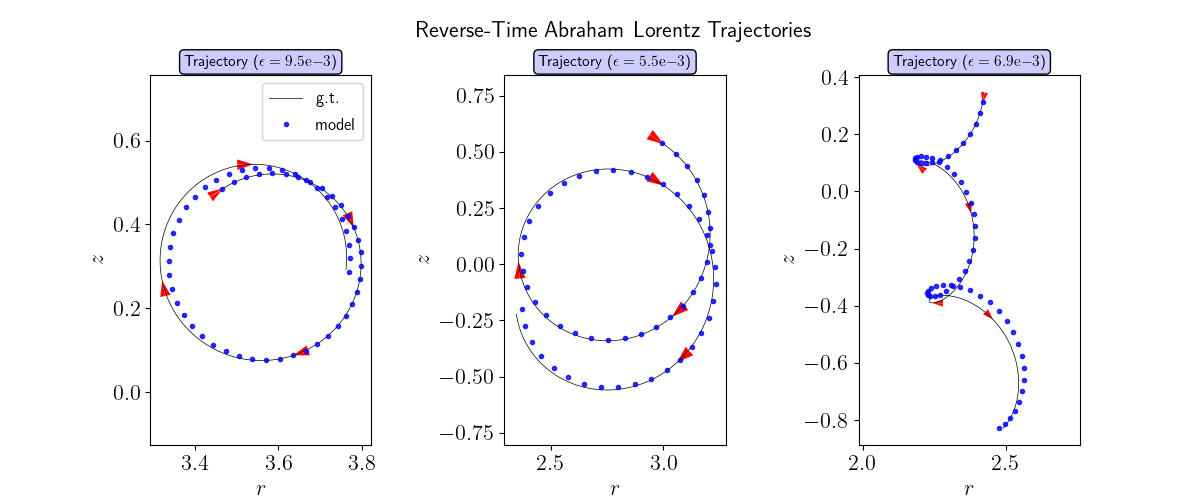}
    \includegraphics[width=\textwidth]{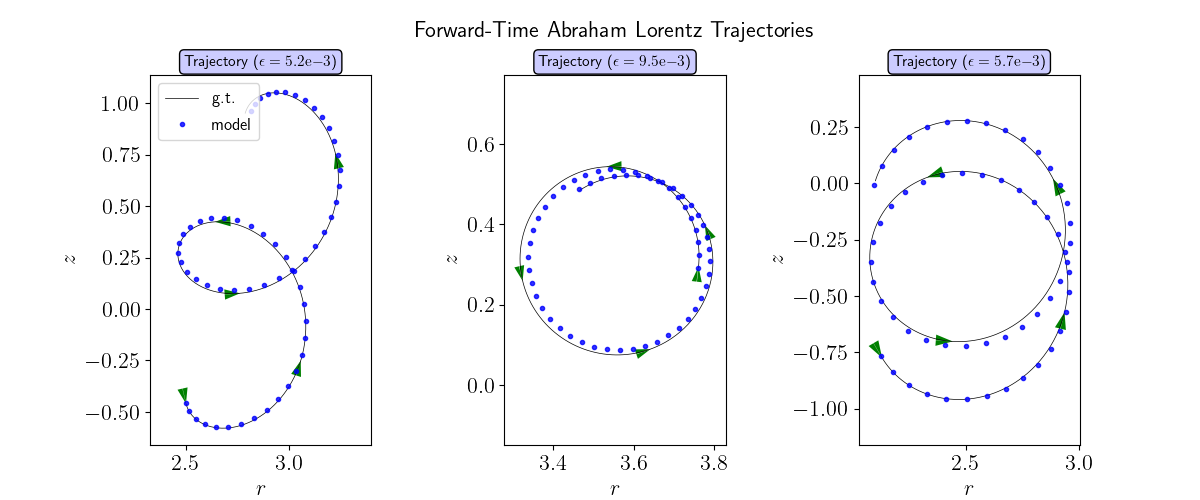}
    \caption{Displayed are the trajectories for three values of $\epsilon$ evolved backwards in time along the slow manifold (top three plots, red arrows) and three values of $\epsilon$ evolved forwards in time along the slow manifold (bottom three plots, green arrows). 
    The solid lines represent the ground truth trajectories and the dotted lines represent the model.} 
    \label{fig:trajectories}
\end{figure}

\begin{figure}[tbh]
    \centering
    \includegraphics[width=\textwidth]{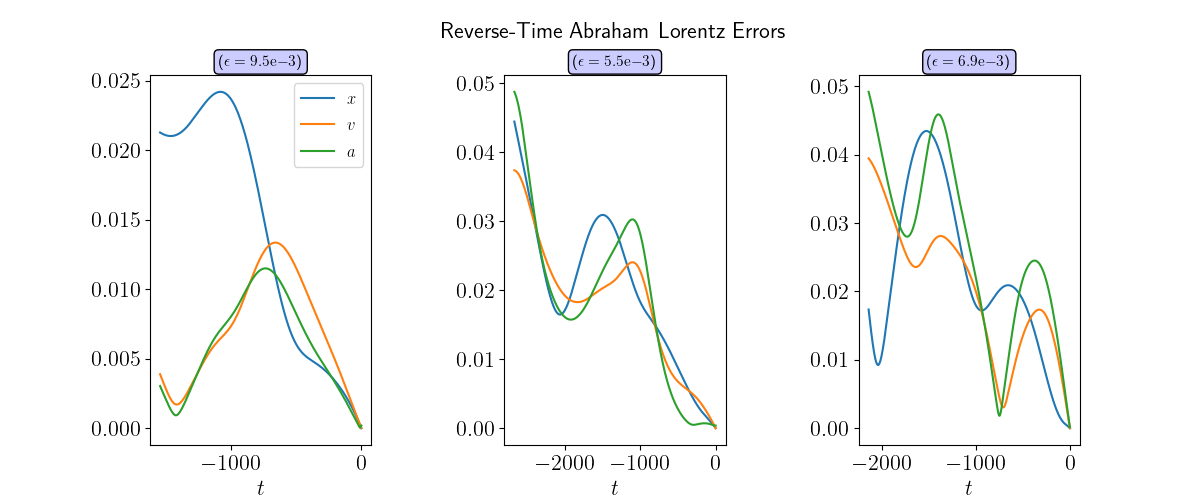}
    \includegraphics[width=\textwidth]{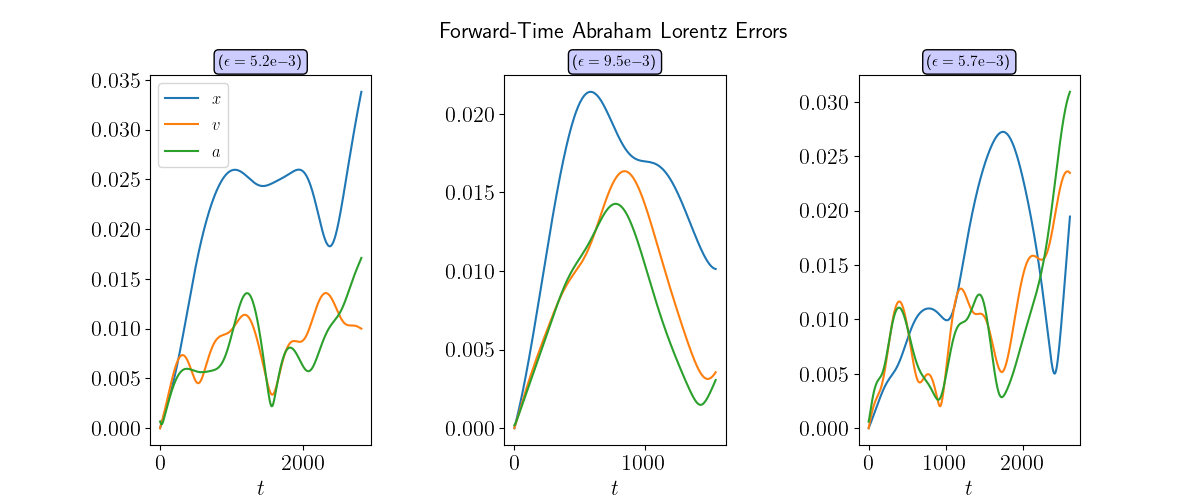}
    \caption{Error trajectories for the reverse time (top three plots) and forward time simulations are displayed. For each simulation, the $L_2$ norm of position, velocity and acceleration are plotted.} 
    \label{fig:errors}
\end{figure}

Figure~\ref{fig:trajectories} shows comparisons 
of the solutions in $r$-$z$ space
between the trained network and trajectories computed using RK45
on the reverse-time system.
Using the reverse-time trajectories, endpoints were 
supplied as initial conditions to the trained network
to predict forward in time. 
For both the reverse-time predictions (red arrows)
and forward-time predictions (green arrows), 
the trained model agrees well with the solution computed 
numerically.
The corresponding $L_2$ errors are shown for 
position, velocity, and acceleration in 
Figure~\ref{fig:errors}, demonstrating that errors
are low when integrating on the slow time scale.

\begin{table}[tbh]
\caption{Comparison of the true eigenvalues {Jacobian of the fast dynamics} to the
eigenvalues of the {negative Schur form} of the trained model for the Abraham-Lorentz problem.\label{tab:eigenvaluesabl}
}
\begin{center}
\begin{tabular}{|>{\centering\arraybackslash}p{2.5cm}|>{\centering\arraybackslash}p{2.5cm}|>{\centering\arraybackslash}p{2.5cm}|}
\hline
\multicolumn{3}{|c|}{Reverse-Time Abraham-Lorentz System Eigenvalues} \\
\hline
 g.t. &   model &      error \\
 \hline
     -1 & -0.99962 &    3.8e-04 \\
     -1 & -0.99958 &    4.2e-04 \\
     -1 & -0.99958 &    4.2e-04 \\
\hline
\end{tabular}
\end{center}
\end{table}

Table~\ref{tab:eigenvaluesabl} shows the comparison of
the true eigenvalues of the linear part of the fast dynamics 
of the reverse-time system. The eigenvalues were recovered
with a great deal of accuracy, demonstrating that the 
fast-scale dynamics are learned accurately.

\begin{figure}[tbh]
    \centering
    \includegraphics[width=\textwidth]{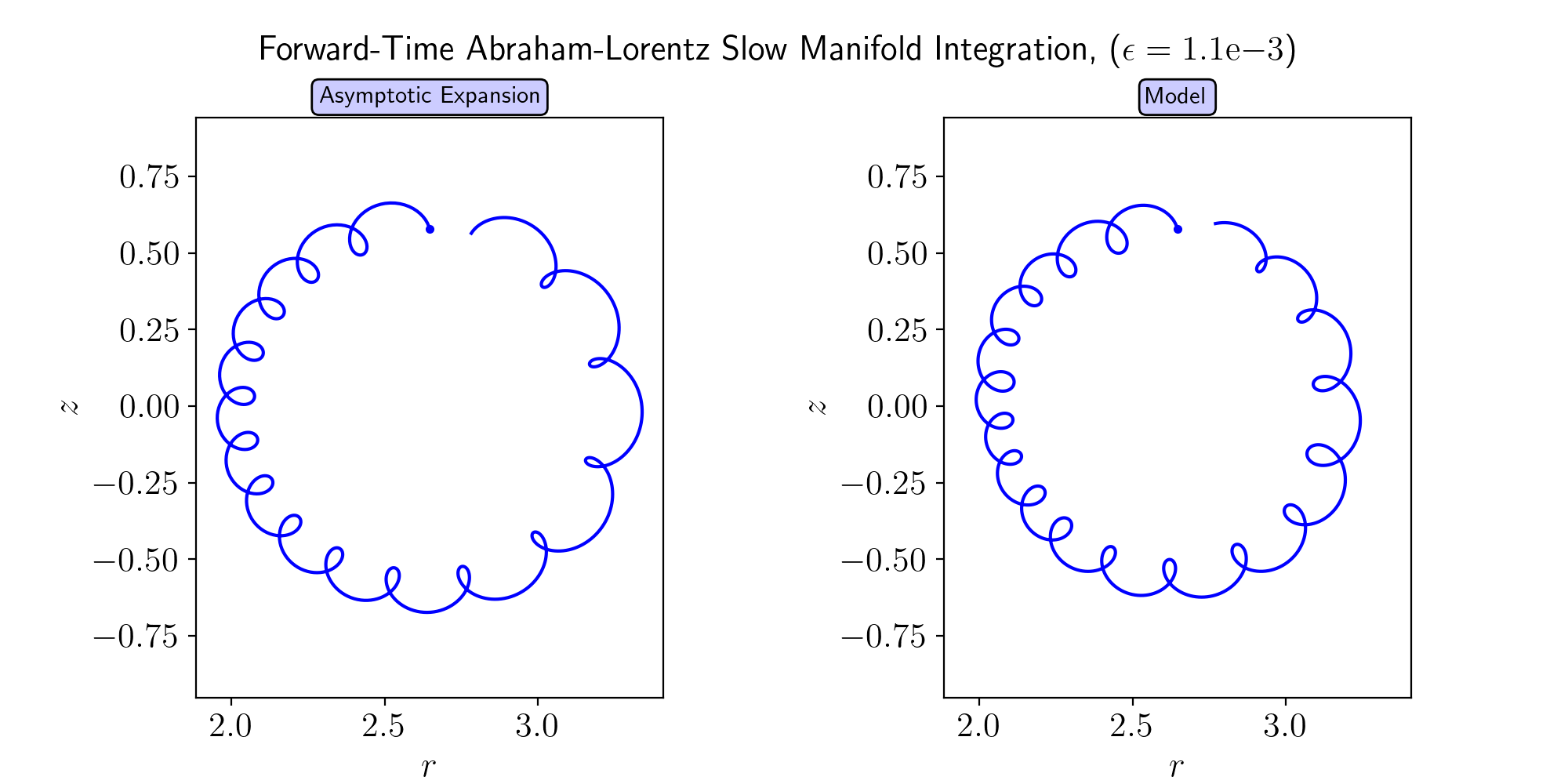}
    \includegraphics[width=\textwidth]{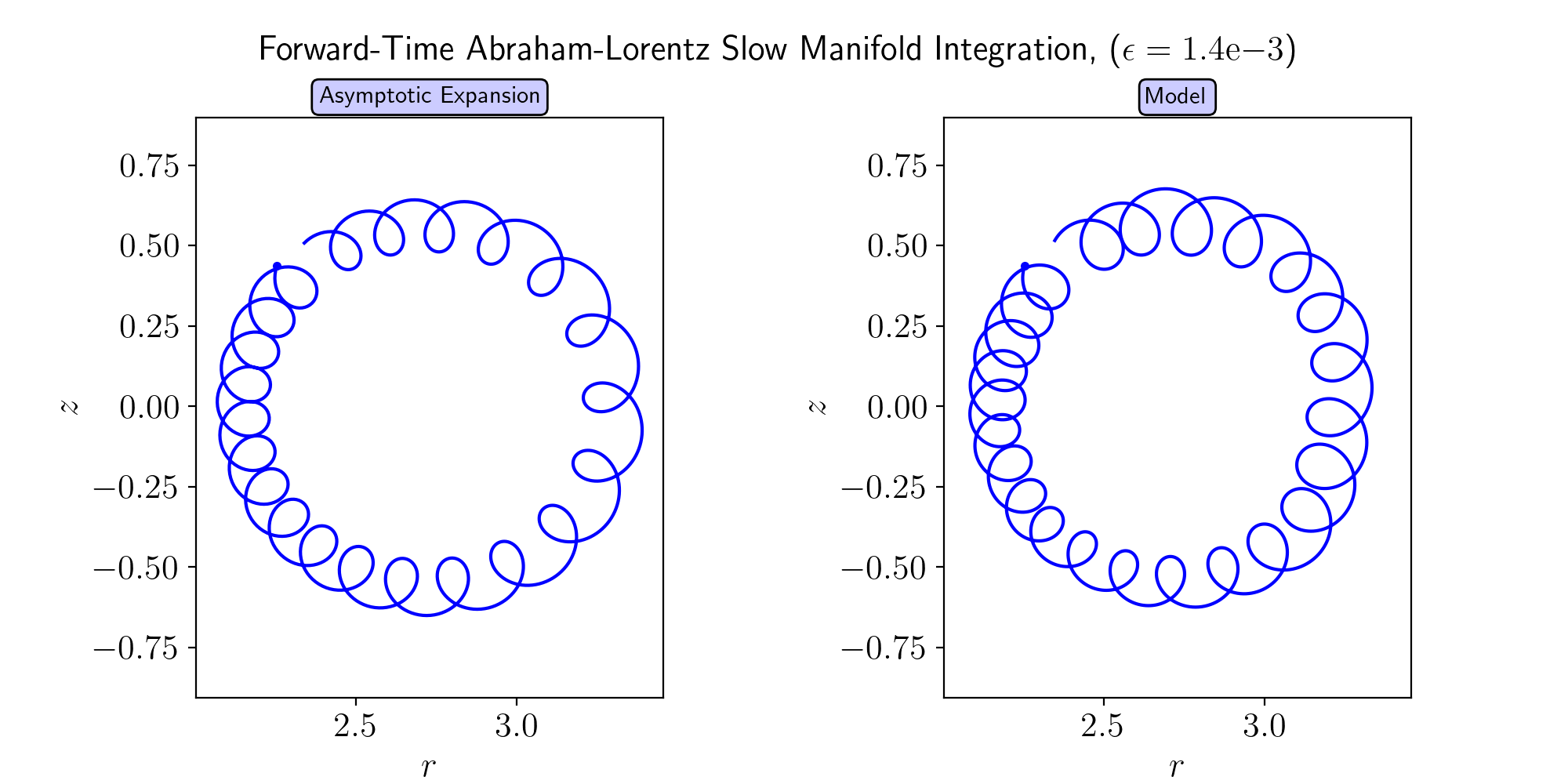}
    \caption{Comparison of orbits of the forward-time Abraham-Lorentz system~\eqref{eq:forward} along the slow manifold for
    the first two terms of the asymptotic expansion 
    in~\eqref{eq:expansion} (left) and the 
    trained model (right) for two different initial 
    conditions and choices of $\epsilon$.}
    \label{fig:orbits}
\end{figure}

Figure~\ref{fig:orbits} shows comparisons of 
orbits between the trained model and orbits obtained 
when integrating the system by substituting 
the first two terms of the asymptotic expansion 
in~\eqref{eq:expansion}. 
Note again that when integrating along the slow manifold, only the slow variable is evolved using the transformation $h$ and the closure $y=0$.
The trained model is able to qualitatively capture
the gyro-motion behavior 
while also accurately capturing the 
guiding center dynamics over the slow time scale.
Despite training using {trajectories
over a relatively short time frame} (2 or less gyration periods)
due to the slow time scale instability of the reverse-time system,
the trained network is highly capable of predicting long-time
dynamics in forward-time for many more gyration cycles.
We remind the reader that standard numerical methods
cannot integrate the forward-time Abraham-Lorentz system
along the slow manifold without experiencing an exponential growth
instability. 
Direct integration of the asymptotic expansion removes the
instability, but introduces an error dependent on $\epsilon$
due to the asymptotic approximation.
This result highlights the unique capability of the proposed network.

\section{Conclusions}
\label{sec:conclusions}

In this paper, we introduced the fast-slow neural network (FSNN),
which is a data-driven approach suitable for learning
singularly perturbed dynamical systems where the fast-scale
dynamics are dissipative.
Our method enforces the existence of a trainable, attracting invariant slow manifold as a hard constraint.
Invertible neural networks and neural ODEs are 
key components of the FSNN. 
We introduce bi-Lipschitz affine transformation (bLAT) layers for enforcing regularity of the invertible neural network. 
We also develop a negative Schur form network which parameterizes matrices with eigenvalues to the left of the imaginary axis. 
This in addition to the use of an additive L-stable diagonally-implicit Runge Kutta scheme for integration
are key pieces in ensuring that the fast-scale dynamics are stable.
The Schur form network also enables the use of back-substitution 
to efficiently handle the implicit solve.
Analytical representation of the slow manifold 
enables efficient integration on the slow time scale.
We demonstrate the FSNN on many examples that exhibit
two timescales, including the Grad moment system from hydrodynamics,
two-scale Lorenz96 equations for modeling atmospheric dynamics,
and Abraham-Lorentz dynamics modeling radiation reaction of electrons in a magnetic field.

{In order to effectively train the FSNN, a sufficient number of trajectories 
across multiple values of $\epsilon$ and initial conditions are necessary.
Despite this, 
we demonstrate that the network can be used to discover a full-order model with long-term accuracy near a slow manifold when only short duration trajectories
are available. 
Therefore the proposed network is also applicable to discover a reduced-order model as a unique closure discovery procedure. }

{In the examples we consider, the equations and slow manifold reduction in the limit
$\epsilon=0$ are known a priori.
Therefore, we were able to bypass an important shortcoming of our approach, which is 
the need to determine the dimensionality of the slow manifold.
While a plethora of techniques exist in the literature~\cite{manifoldlearning}, 
it would be interesting to combine them to our technique to tackle problems
where the dimensionality is unknown.}

{A future direction of our work involves considering hyperbolic fast slow systems,
which can be represented using a more general Fenichel normal form consisting 
of both stable and unstable directions.}

\section*{Acknowledgment}
This research used resources provided by the Los Alamos National Laboratory Institutional Computing Program, which is supported by the U.S.~Department of Energy National Nuclear Security Administration under Contract No. 89233218CNA000001, and the National Energy Research Scientific Computing Center (NERSC), a U.S. Department of Energy Office of Science User Facility located at Lawrence Berkeley National Laboratory, operated under Contract No. DE-AC02-05CH11231 using NERSC award ASCR-ERCAP0023112.

\bibliographystyle{elsart-num}
\bibliography{references}
\end{document}

%% file: defs.tex
\newcommand{\citeCount}[1]{}
\newcommand{\dt}{\Delta t}

\def\ba#1\ea{\begin{align}#1\end{align}}

\def\bas#1\eas{\begin{align*}#1\end{align*}}

\def\bat#1\eat{\begin{alignat}{3}#1\end{alignat}}

\def\bats#1\eats{\begin{alignat*}{3}#1\end{alignat*}}

\newcommand{\bse}{\begin{subequations}}
\newcommand{\ese}{\end{subequations}}

\newcommand{\av}{\mathbf{ a}}

\newcommand{\uv}{\mathbf{ u}}

\newcommand{\xv}{\mathbf{ x}}

\newcommand{\Bv}{\mathbf{ B}}

\newcommand{\grad}{\nabla}

\clearpage

%% file: trimFig.tex
\newlength{\tfwidth}
\newlength{\tfheight}
\newlength{\tfxa}
\newlength{\tfxb}
\newlength{\tfya}
\newlength{\tfyb}
\newcommand{\trimFigWithBox}[6]{%
\setlength\fboxsep{0pt}%
\setlength\fboxrule{1.0pt}
\fbox{\includegraphics[width=#2, clip, trim=#3 #4 #5 #6]{#1}}%
}
\newcommand{\trimFigNoBox}[6]{%
\setlength\fboxsep{1pt}
\setlength\fboxrule{0.0pt}
\fbox{\includegraphics[width=#2, clip, trim=#3 #4 #5 #6]{#1}}%
}
\newcommand{\trimFigHeightWithBox}[6]{%
\setlength\fboxsep{0pt}%
\setlength\fboxrule{1.0pt}
\fbox{\includegraphics[height=#2, clip, trim=#3 #4 #5 #6]{#1}}%
}
\newcommand{\trimFigHeightNoBox}[6]{%
\setlength\fboxsep{1pt}
\setlength\fboxrule{0.0pt}
\fbox{\includegraphics[height=#2, clip, trim=#3 #4 #5 #6]{#1}}%
}

\newsavebox\figBox

\newcommand{\trimw}[6]{%
\sbox\figBox{\includegraphics{#1}}
\setlength{\tfwidth}{\the\wd\figBox}
\setlength{\tfheight}{\the\ht\figBox}
\setlength{\tfxa}{\tfwidth*\real{#3}}%
\setlength{\tfxb}{\tfwidth*\real{#4}}%
\setlength{\tfya}{\tfheight*\real{#5}}%
\setlength{\tfyb}{\tfheight*\real{#6}}%
\trimFigNoBox{#1}{#2}{\tfxa}{\tfya}{\tfxb}{\tfyb}%
}

\newcommand{\trimwb}[6]{%

\sbox\figBox{\includegraphics{#1}}
\setlength{\tfwidth}{\the\wd\figBox}
\setlength{\tfheight}{\the\ht\figBox}
\setlength{\tfxa}{\tfwidth*\real{#3}}%
\setlength{\tfxb}{\tfwidth*\real{#4}}%
\setlength{\tfya}{\tfheight*\real{#5}}%
\setlength{\tfyb}{\tfheight*\real{#6}}%
\trimFigWithBox{#1}{#2}{\tfxa}{\tfya}{\tfxb}{\tfyb}%
}

\newcommand{\trimh}[6]{%
\sbox\figBox{\includegraphics{#1}}
\setlength{\tfwidth}{\the\wd\figBox}
\setlength{\tfheight}{\the\ht\figBox}
\setlength{\tfxa}{\tfwidth*\real{#3}}%
\setlength{\tfxb}{\tfwidth*\real{#4}}%
\setlength{\tfya}{\tfheight*\real{#5}}%
\setlength{\tfyb}{\tfheight*\real{#6}}%
\trimFigHeightNoBox{#1}{#2}{\tfxa}{\tfya}{\tfxb}{\tfyb}%
}

\newcommand{\trimhb}[6]{%

\sbox\figBox{\includegraphics{#1}}
\setlength{\tfwidth}{\the\wd\figBox}
\setlength{\tfheight}{\the\ht\figBox}
\setlength{\tfxa}{\tfwidth*\real{#3}}%
\setlength{\tfxb}{\tfwidth*\real{#4}}%
\setlength{\tfya}{\tfheight*\real{#5}}%
\setlength{\tfyb}{\tfheight*\real{#6}}%
\trimFigHeightWithBox{#1}{#2}{\tfxa}{\tfya}{\tfxb}{\tfyb}%
}